%% file: article.tex
\documentclass[hidelinks,onefignum,onetabnum]{siamart250211}
\input{ex_shared}

\ifpdf
\hypersetup{
  pdftitle={DF-MartNet},
  pdfauthor={Wei Cai, 
    Shuixin Fang, 
    Wenzhong Zhang, 
    Tao Zhou}
}
\fi

\begin{document}

\maketitle

% REQUIRED
\begin{abstract}
    We propose a martingale deep learning method for very high-dimensional stochastic optimal control problems (SOCPs) via their associated Hamilton--Jacobi--Bellman equations and the verification theory for the optimal control. 
    The method decomposes the problem into two coupled components: a martingale formulation for the value function and an integral optimality condition for the feedback control. 
    For the value function, we extend DeepMartNet by introducing a pilot process for state-space exploration and system processes for martingale property enforcement. 
    This yields a derivative-free formulation that avoids computing PDE solution's derivatives and online simulations of full trajectories, while enabling parallel loss evaluation in both time and space. 
    Also, the martingale condition is imposed in a weak Galerkin form realized with adversarial learning, which avoids direct computation of conditional expectations. 
    The control is learned by minimizing an integral residual of the optimality condition over the explored state space, thereby avoiding state-space point-wise optimization over the control set. 
    Numerical results demonstrate accurate and efficient performance for SOCPs with complex patterns in value function  and controls in dimensions up to 10,000.
\end{abstract}
% In this paper, a highly parallel and derivative-free martingale neural network method, based on the probability theory of Varadhan's martingale formulation of PDEs, is proposed to solve stochastic optimal control problems (SOCPs) through their associated Hamilton-Jacobi-Bellman (HJB) equations. 
% The HJB equations are reformulated into a martingale problem such that loss functions will not require the computation of the gradient or Hessian matrix of the value function, and can be computed in parallel in both time and spatial domains. 
% Moreover, the martingale conditions for the equations are enforced using a Galerkin method realized with adversarial learning techniques, eliminating the need for direct computation of the conditional expectations associated with the martingale property. 
% The proposed method is also applicable to quasilinear parabolic PDEs.    

% REQUIRED

% REQUIRED

\section{Introduction}\label{sec_intro}

Stochastic optimal control problems (SOCPs) concern the control of dynamical systems subject to random perturbations, with the goal of optimizing a prescribed performance criterion \cite{Pham2009Continuous}. 
These problems arise widely in decision-making and systems under uncertainty, with applications in engineering, economics, biology, and rare event simulation \cite{Dupuis2012Importance,Vanden2012Rare}; see, e.g., \cite{Zarip2001solution,Nadtochiy2013approximation,Dain2026Uncertainty}. 
Their analysis and numerical solution have long been central topics in applied mathematics.

The theory of SOCPs has evolved along two pioneering and complementary directions: 
Pontryagin's maximum principle \cite{Boltyanskij1961theory,Kushner1965Stochastic,Bismut1978introductory,Bensoussan1983Stochastic,Peng1990general}, which provides necessary conditions for optimal controls via adjoint processes \cite[p.~118]{Yong1999Stochastic}, and Bellman's dynamic programming \cite{Bellman1952On,Bellman1957Dynamic,Bertsekas1976Dynamic,Bellman2010Dynamic}, which characterizes optimal controls through stochastic verification theorems \cite{Zhou1997Stochastic} and reduces the problem to solving the Hamilton--Jacobi--Bellman (HJB) equation \cite[p.~163]{Yong1999Stochastic}.
After decades of development, these two perspectives were connected with the theory of backward stochastic differential equations (BSDEs), 
%a class of stochastic differential equations 
a concept originally introduced by Bismut \cite{Bismut1976Probabilistic} as adjoint equations in the stochastic maximum principle (SMP), and later established as a systematic theory by Pardoux and Peng \cite{pardoux1990adapted}.
Nowadays, the theory of SOCPs continues to be an active area of research, with ongoing advances in analysis, probability, and in more general problem settings \cite{Peng1999Fully,Buckdahn2014Nonlinear,Wang2018Optimal,Carmona2018Probabilistic,Hu2022Global,Ji2024Stochastic}.

Numerical methods for SOCPs play a central role in bridging theoretical developments and practical applications. 
However, designing efficient solvers remains highly challenging in high-dimensional settings, which frequently arise in applications involving large state and control spaces.
It is well-known that conventional numerical methods \cite{Smears2014Discontinuous,osher1991shu,Richardson2006Numerical,Barles2002On,Cacace2012patchy,zhaosweep,zhao2006new,zhao2010stable,zhao2014new,Gong2017efficient,Fu2020Highly,Wang2022Sinc} suffer from the curse of dimensionality (CoD) \cite{Bellman1957Dynamic}: their computational complexity grows exponentially with the problem dimension, making them impractical for high-dimensional SOCPs.

In recent years, the incorporation of deep learning into numerical methods has led to substantial progress in tackling the CoD, yielding practical solvers for high-dimensional SOCPs. 
Existing deep learning methods for SOCPs can be broadly categorized as follows:
(a) SMP approach: the SMP characterizes the optimal control through a coupled forward-backward stochastic differential equations (FBSDEs) system together with a maximum condition \cite[p.~118]{Yong1999Stochastic}, and this system is then solved by deep learning; see, e.g., \cite{cai2024socmartnet,Ji2022Solving,Huang2025Convergence,Ji2025novel}.
(b) PDE approach: the dynamic programming principle reduces the SOCP to solving the associated HJB equation \cite[p.~163]{Yong1999Stochastic}, which is then tackled by deep-learning-based PDE solvers, e.g., \cite{Raissi2019Physics,han2018solving,beck2019deep,wang2022is,Darbon2021some,Darbon2020Overcoming,Al2022Extensions,Ruthotto2020machine}.
(c) Hybrid approach: SMP- and PDE-based approaches are combined to improve state-space exploration and learning efficiency; see, e.g., \cite{Li2024Neural,Onken2023neural}.
(d) Reinforcement learning approach: the SOCP is reformulated as a Markov decision process and solved using policy-iteration methods; see, e.g., \cite{Zhou2021Actor,Yanwei2022Policy,Yanwei2023Policy,Hure2021Deep,Bachouch2022Deep,Zhou2024Solving,Zhou2025Policy}.

Despite the notable success of deep learning approach, existing methods still face substantial challenges in very high-dimensional settings, for example when the state dimension $d$ or control dimension $m$ is in the order of thousands or more.
These challenges arise mainly from the following three sources:
(i) HJB equations involve Hessian matrices with $d \times d$ entries, and PDE- and hybrid-based approaches typically evaluate these Hessians through automatic differentiation of deep neural networks, which become prohibitively costly in computation and memory as the dimension increases. 
This bottleneck has motivated a line of recent works on high-dimensional PDEs \cite{He2023Learning,hu2024sdgd,Hu2024Hutchinson,shi2024stochastic,Li2024computational,Li2024DOF, ZhangHighDimReview2026}.
(ii) SMP- and reinforcement-learning-based approaches avoid explicit Hessian evaluation through reformulation, but they require simulation of sample paths of the state process. 
These simulations are sequential in time and must be repeatedly updated as the control changes during training, which substantially limits the parallelizability of the algorithms.
(iii) Both the SMP and HJB formulations involve an embedded optimization problem for the optimal control, whose minimizer generally does not admit a closed-form expression. 
As a result, pointwise optimization is required at every time-space point, which further exacerbates the CoD problem.

In this work, we develop a new deep learning method for very high-dimensional SOCPs through their associated HJB equations and the
verification theory for the optimal control.
Our method aims to overcome the above-mentioned difficulties by combining the following ideas:
\begin{enumerate}
    \item
    Following the approach of reinforcement learning, we split the HJB equation into two coupled subproblems: one for the value function and one for the optimal control.
    For the value-function subproblem, we extend DeepMartNet \cite{cai2023deepmartnet,cai2023deepmartnet2} to reformulate the equation as a martingale problem, therefore avoiding all derivative computations in the original PDE.

    \item
    To improve the parallelizability of the martingale formulation, we introduce two types of stochastic processes: a pilot process to explore the state space and system processes to encode the martingale property of the value-function subproblem.
    The pilot process is decoupled from training and can be generated offline.
    The system processes can be simulated independently at different time steps, enabling mini-batch sampling and parallel loss evaluation in both time and space.

    \item
    Following the idea of SOC-MartNet \cite{cai2024socmartnet}, we enforce the martingale formulation with a Galerkin method, which is realized with adversarial learning techniques \cite{Zang2020Weak}.
    This avoids computation of conditional expectations in the martingale formulation, and further improves the efficiency of the method.

    \item
    For the optimal-control subproblem, we introduce a neural network to learn the optimal control.
    Its loss is designed to minimize the integrated residual of the optimality condition over the state space explored by the pilot process. 
    This approach eliminates the need for pointwise optimization at every time-space point.
\end{enumerate}
Combining these ideas yields a martingale deep learning method for SOCPs and HJB equations. 
Numerical results show that the proposed method efficiently solves problems, including cases in which the optimal control does not admit an explicit expression, in dimensions up to $10^4$.

The rest of this paper is organized as follows.
Section~\ref{sec_Pre} presents the setup of the SOCP and a brief review of the dynamic programming principle.
Section~\ref{sec_martPDE} develops the martingale deep learning method for SOCPs and the associated HJB equations.
Section~\ref{sec_numres} presents the numerical results. 
Some concluding remarks are given in Section~\ref{sec_conclu}. 

% More details about value function of a given control (not necessarily optimal) and the "HJB" equ.

\section{Problem setup and preliminary}\label{sec_Pre}

\subsection{Problem setup}

Let $(\Omega, \mathcal{F}, \mathbb{P})$ be a complete probability space equipped with a filtration $(\mathcal{F}_t)_{t \in [0, T]}$ generated by the standard Brownian motion $B: [0, T] \times \Omega \to \R^q$, where $T > 0$ is a determined terminal time.
The state equation is given by the following stochastic differential equation (SDE):
\begin{equation}\label{eq_sde}
    X_t = X_0 + \int_0^t \mu\vbr{\big}{s, X_s, u(s, X_s)} \di s + \int_0^t \sigma\vbr{\big}{s, X_s, u(s, X_s)} \di B_s,\quad t \in [0, T],
\end{equation}
where $X_0: \Omega \to \R^d$ is an $\mathcal{F}_0$-measurable random variable with a prescribed distribution; 
the stochastic integral with respect to $B_s$ is understood in the It\^o sense; 
$\mu: [0, T] \times \R^d \times U \to \R^d$ and $\sigma: [0, T] \times \R^d \times U \to \R^{d \times q}$ are given functions representing the drift and diffusion coefficients, respectively;
$U \subset \R^m$ is a given control domain; 
$X: [0, T] \times \Omega \to \R^d$ is called the state process corresponding to the state-feedback control $u: [0, T] \times \R^d \to U$, which belongs to the admissible set $\mathcal{U}_{\mathrm{ad}}$ defined by 
\begin{equation}\label{eq1_defUad}
    \mathcal{U}_{\mathrm{ad}} := \bbr{u: [0, T] \times \R^d \to U: u \text{ is Borel measurable}}.
\end{equation}
Here the Borel measurability of $u$ is required to ensure the validity of the relevant stochastic integrals and expectations in the subsequent analysis. 
To avoid technicalities, we assume that for any $u \in \mathcal{U}_{\mathrm{ad}}$, the SDE~\eqref{eq_sde} admits a unique strong solution $X$. 
Sufficient conditions for this assumption can be found in, e.g., \cite[Chapter~1, Theorem~6.17]{Yong1999Stochastic}.

The performance of a control $u \in \mathcal{U}_{\mathrm{ad}}$ is measured by the cost functional
\begin{equation}\label{eq_costfunc}
    J(u) := \E{\int_0^T c\vbr{\big}{s, X_s, u(s, X_s)} \di s + g(X_T)}, 
\end{equation}
where $c: [0, T] \times \R^d \times U \to \R$ and $g: \R^d \to \R$ are given functions representing the running cost and the terminal cost, respectively.
The SOCP considered by this paper is to find an optimal control $u^* \in \mathcal{U}_{\mathrm{ad}}$ minimizing the cost functional, i.e.,
\begin{equation}\label{eq_SOCP}
    J(u^*) = \inf_{u \in \mathcal{U}_{\mathrm{ad}}} J(u).
\end{equation}

\subsection{Dynamic programming principle}

To lay the theoretical foundation for our method, we briefly review the dynamic programming principle in stochastic control; see, e.g., \cite[Chapter 4]{Yong1999Stochastic} and \cite[Chapter 3]{Pham2009Continuous} for comprehensive treatments.

The dynamic programming principle is based on the value function defined by
\begin{equation}\label{eq_defval_u}
    v(t, x) := \inf_{u \in \mathcal{U}_{\mathrm{ad}}} \E{\int_t^T c\vbr{\big}{s, X_s, u(s, X_s)} \di s + g(X_T) \;\Big\vert\; X_t = x}
\end{equation}
for $(t, x) \in [0, T] \times \R^d$,
where $X$ is the state process governed by \eqref{eq_sde} under control $u$. 
By \cite[p. 182, Proposition 3.5; p. 190, Theorem 5.2]{Yong1999Stochastic}, the value function $v$ is the classical or viscosity solution to the following HJB equation:
\begin{equation}\label{eq_HJBPDE}
\left\{\begin{aligned}
    &\partial_t v(t, x) + \inf_{\kappa \in U} \bbr{\mathcal{L}^{\kappa} v (t, x) + c(t, x, \kappa)} = 0, \quad (t, x) \in [0, T)  \times \R^d,\\
    &v(T, x) = g(x), \quad x \in \R^d,
\end{aligned}\right.
\end{equation}
where $\mathcal{L}^{\kappa}$ is a differential operator defined by
\begin{equation}\label{eq_defLkappa}
    \mathcal{L}^{\kappa} := \mu^{\top}(t, x, \kappa) \partial_x + \frac{1}{2} \operatorname{Tr}[\sigma \sigma^{\top}\br{t, x, \kappa} \partial_{xx}^2], \quad (t, x, \kappa) \in [0, T) \times \R^d \times U.
\end{equation}
The operators $\partial_x := \nabla_x$ and $\partial_{xx}^2 := \nabla_x \nabla_x^{\top}$ denote the gradient and Hessian, respectively, and $\operatorname{Tr}[\cdot]$ denotes the trace of a matrix. 

Under the regularity condition $v \in C^{1,2}$, i.e., $v(t, x)$ is once continuously differentiable in $t \in [0, T)$ and twice continuously differentiable in $x \in \R^d$, the verification theorem \cite[Chapter~5, Theorem~5.1]{Yong1999Stochastic} implies that the optimal control $u^*$ can be given by the following feedback form:
\begin{equation}\label{eq_optcontrol}
    u^*(t, x) = \argmin_{\kappa \in U} \bbr{\mathcal{L}^{\kappa} v (t, x) + c\vbr{\big}{t, x, \kappa}}, \quad (t, x) \in [0, T) \times \R^d.
\end{equation}
Thus, the SOCP~\eqref{eq_SOCP} reduces to solving the HJB equation~\eqref{eq_HJBPDE} for $v$ and then computing $u^*$ from \eqref{eq_optcontrol}.
This paper aims to develop an efficient implementation of this procedure for high-dimensional SOCPs.

\section{Martingale deep learning}\label{sec_martPDE}

To proceed, we substitute \eqref{eq_optcontrol} into \eqref{eq_HJBPDE}, reformulating the HJB equation as the following two-step procedure for determining the value function $v$ and the optimal control $u$:
\begin{align}
    \sptext{find} v \sptext{s.t.} &\bigl(\partial_t + \mathcal{L}^u\bigr) v (t, x) + c\vbr{\big}{t, x, u(t, x)} = 0, \label{eq_pt_Lu} \\
    \sptext{find} u \sptext{s.t.} &u(t, x) = \argmin_{\kappa \in U} \bbr{\mathcal{L}^{\kappa} v (t, x) + c\vbr{\big}{t, x, \kappa}}, \label{eq1_optcontrol}
\end{align}
where $(t, x)$ ranges over $[0, T) \times \R^d$; $v$ satisfies the terminal condition in \eqref{eq_HJBPDE} a priori; $\mathcal{L}^u$ is defined by \eqref{eq_defLkappa} with $\kappa = u(t, x)$.

In the remainder of this section, we develop a deep learning method to solve $(u, v)$ simultaneously from \eqref{eq_pt_Lu} and \eqref{eq1_optcontrol}.
The primary difficulty is to evaluate the derivatives in $\mathcal{L}^{\kappa} v$, which involves the Hessian $\partial_{xx}^2 v$ with $d \times d$ entries.
A second difficulty is the minimization in \eqref{eq1_optcontrol}, which generally lacks a closed-form minimizer; a direct approach requires pointwise optimization at every $(t, x)$.
These two difficulties incur substantial computational costs in high dimensions, constituting a primary source of the CoD.
We address them with a martingale approach that is derivative-free, highly parallel, and free of pointwise minimization.

\subsection{Stochastic processes for martingale formulation}

Let $\pi_N$ be a uniform time partition of $[0, T]$ defined by
\begin{equation}\label{eq_timepart}
    \pi_N: 0 = t_0 < t_1 < \cdots < t_N = T, \quad t_n := n \Delta t, \quad \Delta t := T / N
\end{equation}
with $N$ the number of time steps.
Define the Brownian increments by
\begin{equation*}
    \Delta B_{n+1} := B_{t_{n+1}} - B_{t_n}, \quad n = 0, 1, \cdots, N-1,
\end{equation*}
where the subscript $n+1$ in $\Delta B_{n+1}$ is used to indicate that it is $\mathcal{F}_{t_{n+1}}$-measurable.

Let $u_{\alpha}$ and $v_{\theta}$ denote neural networks parameterized by $\alpha$ and $\theta$, respectively, used to approximate the optimal control $u$ and the value function $v$ as defined in \eqref{eq_pt_Lu} and \eqref{eq1_optcontrol}.
The following two types of stochastic processes are the key ingredients for our method:
\begin{enumerate}
    \item {\bf (Pilot process)} 
    Let $\hat{X}: \pi_N \times \Omega \to \R^d$ be a discrete pilot process (indicated by $\hat{} $ ) used to explore the state space $\R^d$, generated by the following Euler--Maruyama scheme:
    \begin{equation}\label{eq_SDE}
        \hat{X}_{n+1} := \hat{X}_n + \hat{\mu}_n \Delta t + \hat{\sigma}_n \Delta B_{n+1}, \quad n = 0, 1, \cdots, N-1,
    \end{equation}
    where $\hat{\mu}_n$ and $\hat{\sigma}_n$ are the drift and diffusion coefficients evaluated with a pilot control $\hat{u}$ in place of $u$, i.e.,
    \begin{equation}\label{eq_hatmu}
        \hat{\mu}_n := \mu\big(t_n, \hat{X}_n, \hat{u}(t_n, \hat{X}_n)\big), \quad \hat{\sigma}_n := \sigma\big(t_n, \hat{X}_n, \hat{u}(t_n, \hat{X}_n)\big).
    \end{equation}
    For parallelization purposes, the pilot path $(\hat{X}_n)_{n=0}^{N}$ need only be updated occasionally during training, with $\hat{u}$ taken as a frozen snapshot of $u_{\alpha}$.
    
    \item {\bf (System process)} 
    For $n = 0, 1, \cdots, N-1$, define the system process as a random jump from $\hat{X}_n$ to a new location $X_{n+1}$ at time $t_{n+1}$ by
    \begin{equation}\label{eq_branch}
        X_{n+1} := \hat{X}_n + \mu_n \Delta t + \sigma_n \Delta B_{n+1}, \quad n = 0, 1, \cdots, N-1, 
    \end{equation}
    where $\mu_n$ and $\sigma_n$ are the drift and diffusion coefficients evaluated using the current control network $u_{\alpha}$, i.e.,
    \begin{equation}\label{eq_mun}
        \mu_n := \mu\big(t_n, \hat{X}_n, u_{\alpha}(t_n, \hat{X}_n)\big), \quad \sigma_n := \sigma\big(t_n, \hat{X}_n, u_{\alpha}(t_n, \hat{X}_n)\big).
    \end{equation}
    The process $(X_{n+1})_{n=0}^{N-1}$ encodes the dynamics of $\mathcal{L}^{u_{\alpha}}$, which will be used to derive a martingale formulation for \eqref{eq_pt_Lu} in \cref{sec_martvalue}. 
\end{enumerate}

\begin{figure}[t]
    \centering
    \subfloat[Conventional approach]
    {\includegraphics[width=0.45\textwidth]{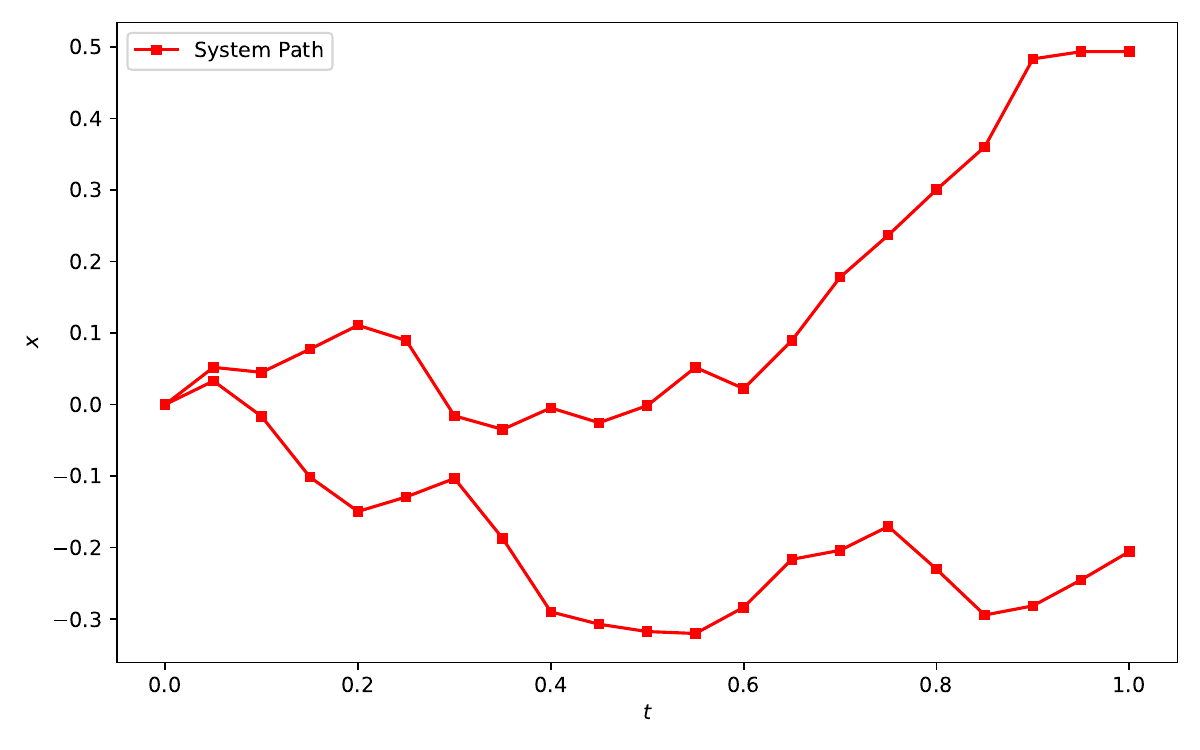}}
    \subfloat[Our approach]{\includegraphics[width=0.45\textwidth]{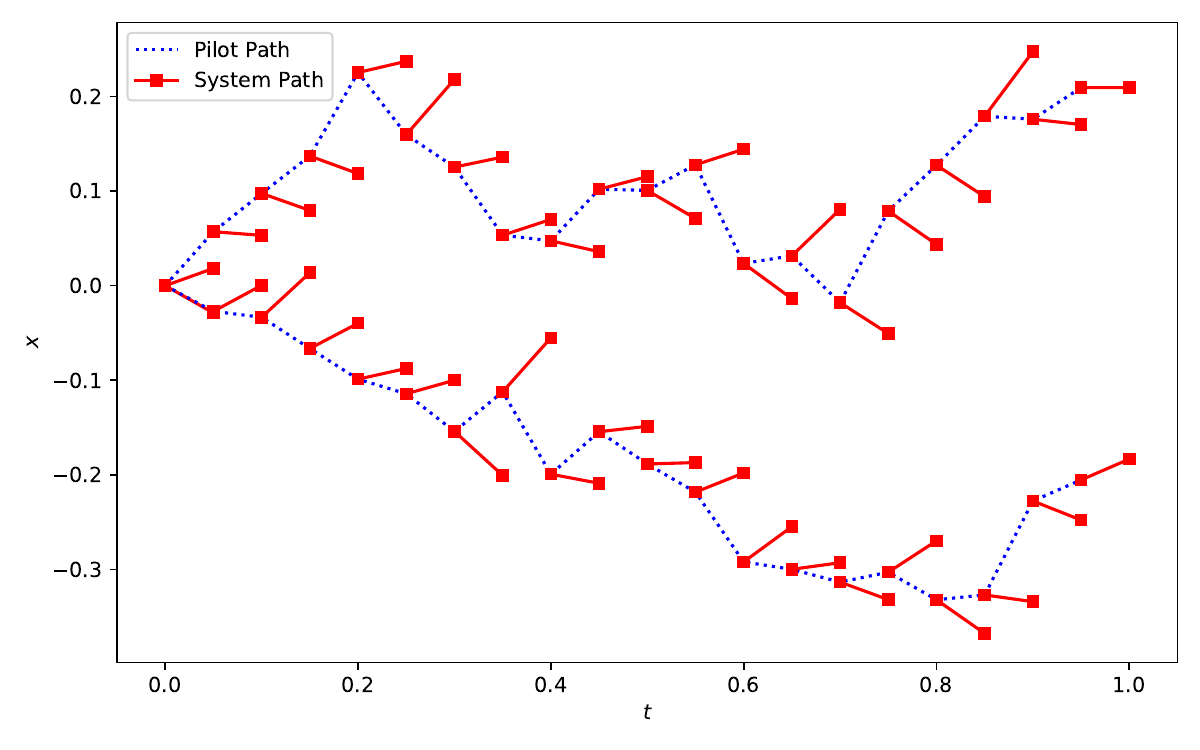}}
    % {\includegraphics[width=0.45\textwidth]{\pathcf sys_paths.pdf}}
    % \subfloat[Our approach]{\includegraphics[width=0.45\textwidth]{\pathcf pilsys_paths.pdf}}
    \caption{
    Comparison of conventional and pilot--system path simulations. The horizontal and vertical axes represent time and state, respectively. In (b), red segments are system-path branches from the pilot path.
    }\label{fig_paths}
\end{figure}

The pilot and system processes have distinct computational features that align with their roles in the method.
The pilot process $\hat{X}$ must be simulated sequentially in time; however, since its main purpose is to explore the state space, it can be generated offline and updated only occasionally during training.
By contrast, the system process $X$ depends on the current control network $u_{\alpha}$ and therefore must be updated at every training iteration to reflect the dynamics of $\mathcal{L}^{u_{\alpha}}$.
Nevertheless, once the pilot path is generated, the system paths $(X_{n+1})_{n=0}^{N-1}$ can be computed independently across all time steps, which enables parallelization in the temporal direction.
These differences are illustrated schematically in \Cref{fig_paths}.

\begin{remark}[Controlled SDEs for space exploration]
    The pilot process $\hat{X}$ in \eqref{eq_SDE} explores the state space by mimicking the controlled dynamics in \eqref{eq_sde}. 
    This strategy helps to focus the training on regions most relevant to the control problem. 
    As shown in \cite{Onken2023neural,Li2024Neural}, such an exploration strategy can significantly improve the performance of SOCP solvers, especially for problems with complex dynamics.
\end{remark}

\subsection{Martingale formulation for value function}\label{sec_martvalue}

We begin by deriving a martingale formulation for the value-function subproblem \eqref{eq_pt_Lu}. 
The main idea is to parallelize the original DeepMartNet \cite{cai2023deepmartnet,cai2023deepmartnet2} using the pilot and system processes introduced above.

The random jump from $\hat{X}_n$ to $X_{n+1}$ in \eqref{eq_branch} can be interpreted as the evolution of the process $\hat{X}_{n, t} := \hat{X}_n + \mu_n (t - t_n) + \sigma_n (B_t - B_{t_n})$ as $t$ increases from $t_n$ to $t_{n+1}$. 
This evolution admits an infinitesimal generator $\partial_t + \mathcal{L}^{u_{\alpha}}$, where the coefficients in $\mathcal{L}^{u_{\alpha}}$ are evaluated at $(t, x) = (t_n, \hat{X}_n)$ and held fixed during the evolution from $t_n$ to $t_{n+1}$.
Then the Dynkin formula \cite[Theorem 7.4.1]{Oksendal2003Stochastic} implies that
\begin{equation}\label{eq_dynkin}
\begin{aligned}
    &\mathbb{E}[v_{\theta}(t_{n+1}, X_{n+1}) \,\big\vert\, \hat{X}_n] - v_{\theta}(t_n, \hat{X}_n) \\
    =\;& \mathbb{E}\Big[ \int_{t_n}^{t_{n+1}} \big(\partial_t + \mathcal{L}^{u_{\alpha}}\big\vert_{(t, x)=(t_n, \hat{X}_n)}\big) v_{\theta}(s, \hat{X}_{n, s}) \di s \,\Big\vert\, \hat{X}_n \Big]\\[0.5ex]
    =\;& (\partial_t + \mathcal{L}^{u_{\alpha}}) v_{\theta}(t_n, \hat{X}_n) \Delta t + \mO(\Delta t^2),
\end{aligned}
\end{equation}
where the remainder term is valid under the regularity condition $v_{\theta} \in C^{1,2}$.

To relate the above formula to the HJB equation \eqref{eq_pt_Lu}, we add the running cost term $c(t, x, u_{\alpha}(t, x)) \Delta t$, evaluated at $(t, x) = (t_n, \hat{X}_n)$, to both sides, which yields
\begin{equation}\label{eq_martResu}
\begin{aligned}
    &\mathbb{E}\big[\M(t_n, \hat{X}_n, X_{n+1}; u_{\alpha}, v_{\theta}) \,\big\vert\, \hat{X}_n\big]\\
    =\;& \big\{(\partial_t + \mathcal{L}^{u_{\alpha}}) v_{\theta} (t_n, \hat{X}_n) + c\vbr{\big}{t_n, \hat{X}_n, u_{\alpha}(t_n, \hat{X}_n)} \big\} \Delta t + \mO(\Delta t^2), 
\end{aligned}
\end{equation}
where $\M$ is a deterministic function defined by
\begin{equation}\label{eq_defM}
    \M\br{t, x, x'; u, v} := v\br{t + \Delta t, x'} - v(t, x) + c(t, x, u(t, x)) \Delta t
\end{equation}
for $(t, x, x') \in [0, T] \times \R^d \times \R^d$.
Neglecting the remainder term $\mO(\Delta t^2)$, the second line represents exactly the residual of \eqref{eq_pt_Lu} evaluated at $(t_n, \hat{X}_n)$. 
Therefore, the term in the first line should vanish, which leads to the following martingale formulation for \eqref{eq_pt_Lu}: train $v_{\theta}$ to satisfy
\begin{equation}\label{eq_martPDE}
    \mathbb{E}\big[\M(t_n, \hat{X}_n, X_{n+1}; u_{\alpha}, v_{\theta}) \,\big\vert\, \hat{X}_n\big] = 0, \quad n = 0, 1, \cdots, N-1,
\end{equation}
where $X_{n+1}$ is defined in \eqref{eq_branch}.

The formulation in \eqref{eq_martPDE} is an alternative characterization of the PDE~\eqref{eq_pt_Lu}, which does not involve any derivatives of $v_{\theta}$.
However, a direct enforcement of \eqref{eq_martPDE} with $v_{\theta}$ requires a two-level Monte Carlo procedure: 
first sampling the pilot path $\hat{X}$, and then, for each realized $\hat{X}_n$, approximating the conditional expectation in \eqref{eq_martPDE} by generating multiple samples of $X_{n+1}$ conditioned on $\hat{X}_n$.
This nested Monte Carlo approach is computationally intensive.

To address this issue, we adopt the idea in \cite{Zang2020Weak,cai2024socmartnet} and enforce \eqref{eq_martPDE} in a weak sense using a Galerkin method.
Let $\rho_{\eta}: [0, T] \times \R^d \to \R^r$ be a neural network parameterized by $\eta$, which serves as a test function in the Galerkin method.
Then \eqref{eq_martPDE} can be enforced in the weak sense by requiring $v_{\theta}$ to satisfy
\begin{equation*}
    \sum_{n=0}^{N-1} \E{\M(t_n, \hat{X}_n, X_{n+1}; u_{\alpha}, v_{\theta}) \rho_{\eta}(t_n, \hat{X}_n)} = 0 \sptext{for any} \eta,
\end{equation*}
where the probability measure in the expectation serves as the weight for the Galerkin method.
The above condition can be equivalently written as the following minimax problem:
\begin{equation}\label{eq_minimax}
    \min_{\theta} \max_{\eta} \big|L(u_{\alpha}, v_{\theta}, \rho_{\eta})\big|^2,
\end{equation}
where the loss function $L$ is defined by
\begin{equation}\label{eq_defLoss}
    L(u_{\alpha}, v_{\theta}, \rho_{\eta}) := \sum_{n=0}^{N-1} \E{\M(t_n, \hat{X}_n, X_{n+1}; u_{\alpha}, v_{\theta}) \rho_{\eta}(t_n, \hat{X}_n)}. 
\end{equation}
Here $\M$ is given by \eqref{eq_defM}; $\hat{X}_n$ and $X_{n+1}$ given by \eqref{eq_SDE} and \eqref{eq_branch}, respectively.
The minimax problem in \eqref{eq_minimax} can be solved via the adversarial training procedure proposed in \cite{Zang2020Weak}; see \cref{sec_learning} for details.

\begin{remark}[Efficiency gain from weak formulation]
    Compared with enforcing \eqref{eq_martPDE} directly, which would require a two-level Monte Carlo procedure, 
    the weak formulation \eqref{eq_minimax} reduces the loss evaluation in \eqref{eq_defLoss} to a single-level procedure.
    Specifically, since the expectation in \eqref{eq_defLoss} is taken with respect to the joint distribution of $(\hat{X}_n, X_{n+1})$, it can be approximated by generating independent identically distributed (i.i.d.) pairs of $(\hat{X}_n, X_{n+1})$, instead of generating multiple samples of $X_{n+1}$ for each realized $\hat{X}_n$.
    This substantially lowers the computational cost.
    A mini-batch version of \eqref{eq_defLoss} in such a single-level manner is given in \cref{eq_defL} for practical implementation.
\end{remark}

\begin{remark}[Trick of vector-valued test function]
    Recall that $\rho_{\eta}$ is constructed as an $\R^r$-valued test function. 
    Accordingly, the loss $L$ defined in \eqref{eq_defLoss} is also $\R^r$-valued, and $|\cdot|$ in \eqref{eq_minimax} denotes the Euclidean norm on $\R^r$.
    Compared with the scalar-valued case $r = 1$, choosing vector-valued test functions, e.g., $r \geq 600$, substantially improves the stability of the adversarial training for \eqref{eq_minimax}, as evidenced in the numerical results in \cite{cai2024socmartnet,Cai2026Deep}.
    A specific structure of $\rho_{\eta}$ is given in \eqref{eq_defrhonet}. 
\end{remark}

\begin{remark}[Connection to derivative-free methods for PDEs]
    The equality in \eqref{eq_dynkin} is a key step of the proposed method, as it eliminates the differential operator $\mathcal{L}^{u_{\alpha}}$ through a martingale differencing argument.
    Similar ideas have also been explored in \cite{Han2020derivative,zhang2025shotgun,Cai2026Deep} for solving high-dimensional parabolic PDEs.
\end{remark}

\subsection{Integral formulation for feedback control}

We next derive an integral form of \eqref{eq1_optcontrol} for learning the optimal feedback control $u$.
The main idea is to use \eqref{eq_martResu} to remove the differential operator $\mathcal{L}^u$ and transform the pointwise minimization in \eqref{eq1_optcontrol} into an integral minimization over the state explored by $\hat{X}$.

Note that \eqref{eq1_optcontrol} can be equivalently written as
\begin{equation}\label{eq_argminLv}
    u(t, x) = \argmin_{\kappa \in U} \bbr{(\partial_t + \mathcal{L}^{\kappa}) v (t, x) + c\br{t, x, \kappa}},
\end{equation}
since the term $\partial_t v(t, x)$ is independent of $\kappa$ and thus does not affect the minimizer
; that is, we seek an optimal feedback control $u$ under the current value function $v$.
With $P_{\hat{X}_n}(\cdot)$ denoting the probability density function of $\hat{X}_n$, the minimization condition~\eqref{eq_argminLv} implies
\begin{equation}\label{eq_argminLv4}
\begin{aligned}
    u =\;& \argmin_{\bar{u} \in \mathcal{U}_{\mathrm{ad}}} \sum_{n=0}^{N-1} \int_{\R^d} \bbr{(\partial_t + \mathcal{L}^{\bar{u}}) v (t_n, x) + c(t_n, x, \bar{u}(t_n, x))} P_{\hat{X}_n}(x) \Delta t  \di x\\
    =\;& \argmin_{\bar{u} \in \mathcal{U}_{\mathrm{ad}}} \sum_{n=0}^{N-1} \E{(\partial_t + \mathcal{L}^{\bar{u}}) v (t_n, \hat{X}_n) + c(t_n, \hat{X}_n, \bar{u}(t_n, \hat{X}_n))} \Delta t,
\end{aligned}
\end{equation}
where the second line is obtained by expressing the integral in the first line as an expectation with respect to $\hat{X}_n$.

By the second line of \eqref{eq_argminLv4} with $(u, v)$ replaced by the neural networks $(u_{\alpha}, v_{\theta})$, 
and then using \eqref{eq_martResu} to express the residual term in the second line of \eqref{eq_argminLv4}, 
we obtain the following integral formulation for learning the optimal feedback control:
\begin{equation}\label{eq_argminLv5}
    u_{\alpha} = \argmin_{\bar{u} \in \mathcal{U}_{\mathrm{ad}}} \sum_{n=0}^{N-1} \E{\M(t_n, \hat{X}_n, X_{n+1}; \bar{u}, v_{\theta})}, 
\end{equation}
where we have neglected the remainder term $\mO(\Delta t^2)$ in \eqref{eq_martResu}.
Using the loss function $L$ defined in \eqref{eq_defLoss}, and constructing $u_{\alpha}$ so that its values lie in $U$ as described in \eqref{eq_utx}, the minimization problem~\eqref{eq_argminLv5} can be solved by training $u_{\alpha}$ via
\begin{equation}\label{eq_argminLv6}
    \min_{\alpha} L(u_{\alpha}, v_{\theta}, 1),
\end{equation}
where $L(u_{\alpha}, v_{\theta}, 1)$ is given by \eqref{eq_defLoss} with $\rho_{\eta}$ replaced by the constant function equal to $1$.
As shown in \eqref{eq_argminLv6}, the integral minimization in \eqref{eq_argminLv5} eliminates the time-state space point-wise minimization in \eqref{eq_argminLv} and reduces the problem to a standard neural network training task. 

By alternating the training of $v_{\theta}$ via \eqref{eq_minimax} and $u_{\alpha}$ via \eqref{eq_argminLv6}, we obtain a unified martingale formulation for solving the SOCP~\eqref{eq_SOCP}.
This formulation is simple, since both subproblems are based on the same loss function $L$ in \eqref{eq_defLoss}, differing only in the choice of test function.
This unified structure simplifies implementation and improves computational efficiency by allowing intermediate quantities to be shared across the two subproblems.

\begin{remark}[Sufficiency of \eqref{eq_argminLv6} for optimal feedback control]
    While \eqref{eq_argminLv4} is derived as a necessary condition for \eqref{eq_argminLv}, it also implies the pointwise optimality in sense of
    \begin{equation*}
        u(t_n, \hat{X}_n) = \argmin_{\kappa \in U} \big\{(\partial_t + \mathcal{L}^{\kappa}) v (t_n, \hat{X}_n) + c(t_n, \hat{X}_n, \kappa)\big\}, \quad n = 0, 1, \cdots, N-1,
    \end{equation*}
    almost surely. 
    This guarantees that $u$ acts as the optimal feedback control in the time-space regions explored by $\hat{X}$ with positive probability. 
    Consequently, by neglecting the $\mO(\Delta t^2)$ remainder in \eqref{eq_martResu}, 
    the minimization problem in \eqref{eq_argminLv6} serves as a sufficient approximation of \eqref{eq_argminLv4} for learning the optimal control over the state space explored by $\hat{X}$, provided $\Delta t$ is sufficiently small.
\end{remark}

\subsection{Network architecture}\label{sec_learning}

Now we specify the architectures of the networks $u_{\alpha}$, $v_{\theta}$, and $\rho_{\eta}$ used in the training procedures \eqref{eq_minimax} and \eqref{eq_argminLv6}.
To avoid technical complications, we assume that $u_{\alpha}$ is constructed such that its output always lies within the control domain $U$.
This can be achieved by applying an appropriate output activation function to enforce the range, or, if $U$ is convex, by projecting the output of a neural network onto $U$:
\begin{equation}\label{eq_utx}
    u_{\alpha}(t, x) = \argmin_{y \in U} \big|\psi_{\alpha}(t, x) - y\big|, \quad (t, x) \in [0, T) \times \R^d,
\end{equation}
where $\psi_{\alpha}$ is a neural network with parameters $\alpha$. 
For regular $U$, such as box constraints of the form $U = \prod_{i=1}^m [a_i, b_i]$ with $-\infty < a_i \leq b_i < +\infty$, this projection can be computed efficiently.

For the network $v_{\theta}$, we assume that the terminal condition in \eqref{eq_HJBPDE} is satisfied by construction. 
This can be achieved by embedding the terminal condition $v(T, x) = g(x)$ directly into the architecture of $v_{\theta}$, e.g., 
\begin{equation}\label{eq_vtx}
    v_{\theta}(t, x) = 
    \begin{cases}
        g(x), & \text{if } t = T, \\
        \phi_{\theta}(t, x), & \text{if } t < T
    \end{cases}
\end{equation}
for $(t, x) \in [0, T] \times \R^d$, where $\phi_{\theta}: [0, T) \times \R^d \to \R$ is a neural network parameterized by $\theta$.
The discontinuity of $v_{\theta}$ at $t = T$ does not cause any issues in training, 
since the network is only evaluated at the discrete time steps in $\pi_N$.

The adversarial network $\rho_{\eta}$ in \eqref{eq_defLoss} plays the role of test functions in classical Galerkin methods for solving PDEs.
By our experimental results, $\rho_{\eta}$ is not necessarily very deep; instead, it can be a shallow network with enough output dimensionality.
Following the multiscale neural network ideas in \cite{Liu2020Multi},
we consider a typical $\rho_{\eta}$ given by
\begin{equation}\label{eq_defrhonet}
    \rho_{\eta}(t, x) = \sin \br{\Lambda \br{W_1 t + W_2 x + b }} \in \R^r, 
\end{equation}
where $\eta := (W_1, W_2, b) \in \R^r \times \R^{r\times d} \times \R^r$ is the trainable parameter; $\Lambda(\cdot)$ is a scale layer defined by 
\begin{equation}\label{eq_defLamb}
    \Lambda(y_1, y_2, \cdots, y_r) = \br{c_1 y_1, c_2 y_2, \cdots, c_r y_r}^{\top} \in \R^r, \quad c_i := 1 + \frac{i-1}{r-1} (c - 1)
\end{equation}
for $y_i \in \R$, where $c > 0$ is a non-trainable hyperparameter; 
$\sin(\cdot)$ is the activation function applied to the outputs of $\Lambda$ in an element-wise manner.

\subsection{Learning algorithm}

The learning algorithm is obtained by approximating the objectives in \eqref{eq_minimax} and \eqref{eq_argminLv6} with mini-batch estimates of the expectation in \eqref{eq_defLoss}.
Specifically, for each $n = 0, 1, \cdots, N-1$, let $\{(\hat{X}_n^m, X_{n+1}^m)\}_{m=1}^M$ be i.i.d.\ samples of $(\hat{X}_n, X_{n+1})$ generated by \eqref{eq_SDE} and \eqref{eq_branch}.  
Then, by the definition of $L$ in \eqref{eq_defLoss}, the objectives in \eqref{eq_minimax} and \eqref{eq_argminLv6} can be approximated by the following unbiased mini-batch estimators:
\begin{align}
    \abs{L(u_{\alpha}, v_{\theta}, \rho_{\eta})}^2 \approx\;& L^{\top}(u_{\alpha}, v_{\theta}, \rho_{\eta}; \mathbb{A}_1) L(u_{\alpha}, v_{\theta}, \rho_{\eta}; \mathbb{A}_2), \label{eq_approxG2}\\[1.2ex]
    L(u_{\alpha}, v_{\theta}, \rho_{\eta}) \approx\;& L(u_{\alpha}, v_{\theta}, \rho_{\eta}; \mathbb{A}_1 \cup \mathbb{A}_2),\label{eq_approxG3}\\
    L(u_{\alpha}, v_{\theta}, \rho_{\eta}; \mathbb{A}_i) :=\;& \frac{1}{\abs{\mathbb{A}_i}} \sum_{m \in \mathbb{A}_i} \sum_{n=0}^{N-1} \M(t_n, \hat{X}_n^m, X_{n+1}^m; u_{\alpha}, v_{\theta}) \rho_{\eta}(t_n, \hat{X}_n^m), \label{eq_defL}
\end{align}
where $\abs{\mathbb{A}_i}$ denotes the cardinality of 
the index sets $\mathbb{A}_i$ for $i=1,2$, sampled from $\{1, 2, \cdots, M\}$, such that $\mathbb{A}_1 \cap \mathbb{A}_2 = \emptyset$.
The disjointness of $\mathbb{A}_1$ and $\mathbb{A}_2$ ensures that the mini-batch estimator in \eqref{eq_approxG2} is unbiased.

In summary, to solve $(u, v)$ from \eqref{eq_pt_Lu} and \eqref{eq1_optcontrol}, we train $(u_{\alpha}, v_{\theta})$ by enforcing \eqref{eq_minimax} and \eqref{eq_argminLv6}.
Specifically, we update $v_{\theta}$ and $u_{\alpha}$ by gradient descent to minimize the objectives in \eqref{eq_minimax} and \eqref{eq_argminLv6}, respectively, and update the test network $\rho_{\eta}$ by gradient ascent to maximize \eqref{eq_minimax}.
In implementation, the objectives and their stochastic gradients are approximated using the mini-batch estimators in \eqref{eq_approxG2}--\eqref{eq_approxG3}.
For efficiency, the pilot paths $\hat{X}_n^m$ are updated only occasionally during training.
\Cref{alg_amnet} summarizes the complete training procedure.

\begin{remark}[Extension to quasilinear parabolic PDEs]
    Observe that the subproblem~\eqref{eq_pt_Lu} reduces to a quasilinear parabolic PDE of $v$ when the control $u$ is replaced by the solution $v$ itself; namely,
    \begin{equation}\label{eq_pde}
        \bigl(\partial_t + \mathcal{L}^v\bigr) v(t, x) + c\vbr{\big}{t, x, v(t, x)} = 0, \quad (t, x) \in [0, T) \times \R^d,
    \end{equation}
    with terminal condition the same as in \eqref{eq_HJBPDE}.
    This observation shows that \Cref{alg_amnet} extends naturally to the quasilinear parabolic PDE~\eqref{eq_pde}. 
    In practice, one removes Line~\ref{lin_alp} and replaces each remaining occurrence of $u_{\alpha}$ with $v_{\theta}$. 
    The resulting output $v_{\theta}$ then provides an approximation of the solution $v$ to \eqref{eq_pde}.
\end{remark}

\begin{algorithm}[t]
    \caption{Martingale deep learning for the SOCP~\eqref{eq_SOCP}}\label{alg_amnet}
    \begin{algorithmic}[1]
        \Require 
        $v_{\theta}$/$u_{\alpha}$/$\rho_{\eta}$: neural networks parameterized by $\theta$/$\alpha$/$\eta$;
        $\delta_0$/$\delta_{1}$/$\delta_{2}$: learning rates for the networks $v_{\theta}$/$u_{\alpha}$/$\rho_{\eta}$;
        $I$: maximum number of iterations;
        $J$/$K$: number of $(\theta, \alpha)$/$\eta$ updates per iteration;
        $I_0$ and $r\%$: interval for updating spatial samples and the percentage of samples updated;
        $M$: total number of sample paths retained in memory. 
        
        \State Initialize the networks $u_{\alpha}$, $v_{\theta}$ and $\rho_{\eta}$
        \State Generate the pilot paths $\{\hat{X}_n^m\}_{n=0}^{N}$, $1 \leq m \leq M$, by \eqref{eq_SDE} with $\hat{u} = u_{\alpha}$. 
        \For{$i = 0, 1, \cdots, I-1$}
        \State Sample disjoint index subsets $\mathbb{A}_1, \mathbb{A}_2$ from $\{1, 2, \cdots, M\}$.
        \State Generate the system paths $\{X_{n+1}^m\}_{n=0}^{N-1}$ for $m \in \mathbb{A}_1 \cup \mathbb{A}_2$ by \eqref{eq_branch}.
        
        \For{$j = 0, 1, \cdots ,J-1$}
        \Comment{$L$ is given by \eqref{eq_defL}}
        \State $\theta \leftarrow \theta - \delta_{0} \nabla_{\theta} \bbr{L^{\top}(u_{\alpha}, v_{\theta}, \rho_{\eta}; \mathbb{A}_1) L(u_{\alpha}, v_{\theta}, \rho_{\eta}; \mathbb{A}_2)}$ 
        \State \label{lin_alp}$\alpha \leftarrow \alpha - \delta_{1} \nabla_{\alpha} L(u_{\alpha}, v_{\theta}, 1; \mathbb{A}_1 \cup \mathbb{A}_2)$
        \EndFor
        \For{$k = 0, 1, \cdots, K-1$}
        \State $\eta \leftarrow \eta + \delta_{2} \nabla_{\eta} \bbr{L^{\top}(u_{\alpha}, v_{\theta}, \rho_{\eta}; \mathbb{A}_1) L(u_{\alpha}, v_{\theta}, \rho_{\eta}; \mathbb{A}_2)}$
        \EndFor
        
        \If{$i$ is divisible by $I_0$} 
        \State Update $r\%$ of the pilot paths $\{\hat{X}_n^m\}_{n=0}^N$ by \eqref{eq_SDE} with $\hat{u}=u_{\alpha}$
        \EndIf
        
        \EndFor
        \Ensure $(v_{\theta}, u_{\alpha})$ as approximations of $(v, u)$. 
    \end{algorithmic}
\end{algorithm}

\subsection{Discussions on advantages}

The martingale deep learning method presented in \Cref{alg_amnet} builds upon our earlier work on DeepMartNet \cite{cai2023deepmartnet,cai2023deepmartnet2} and SOC-MartNet \cite{cai2024socmartnet}, developed for solving PDEs and SOCPs, respectively.  
The differences between the proposed method and existing works are summarized as follows:
\begin{itemize}
    \item 
    Compared with DeepMartNet \cite{cai2023deepmartnet,cai2023deepmartnet2} and 
    other SDE-based methods such as \cite{Ji2022Solving,Huang2025Convergence,Ji2025novel,Zhou2021Actor,Yanwei2022Policy,Yanwei2023Policy,Hure2021Deep,Bachouch2022Deep,Zhou2024Solving,Zhou2025Policy}, the proposed method avoids simulating full trajectories as explained in \Cref{fig_paths}(a). 
    In those existing methods, path simulation is sequential in time. 
    Moreover, because the simulated paths depend on the control neural networks being trained, they must be regenerated whenever the parameters of the networks are updated.  
    This online and temporally sequential path simulation substantially limits the computational efficiency, thus limiting the capability of these methods for very high-dimensional problems.
    
    \item 
    Compared with the original SOC-MartNet \cite{cai2024socmartnet}, the proposed method features a totally derivative-free formulation.
    Although SOC-MartNet already avoids online path simulation, it still depends on automatic differentiation to evaluate derivatives of $v_{\theta}$.
    The proposed method in this paper introduces the random jumps $\hat{X}_n \mapsto X_{n+1}$ in \eqref{eq_branch} to represent the dynamics of $\mathcal{L}^{u_{\alpha}} v_{\theta}$, thereby eliminating the need for explicit derivative evaluations.
    
    \item 
    Compared with PDE solvers, e.g., \cite{Raissi2019Physics,beck2019deep,wang2022is,weinan2017deep,han2018solving,han2018convergence,raissi2018forwardbackward,Zhang2022FBSDE,Nusken2023Interpolating,Han2020Solving}, applied directly to HJB equations,
    our method avoids time and state-space pointwise minimization over the control variable when learning the optimal feedback control. 
    Instead, the control is learned via the integral formulation in \eqref{eq_argminLv5}, similar to standard neural network training by loss minimization. This avoids the complex optimization steps required in SMP and HJB formulations.
\end{itemize}

In summary, the proposed method is derivative-free, avoids sequential path simulation, and eliminates time and state-space  pointwise minimization.
These features yield substantial improvements in computational efficiency, especially for very high-dimensional problems, as demonstrated in \Cref{sec_numres}.

\section{Numerical examples}\label{sec_numres}

We will demonstrate the effectiveness of the proposed method by applying \Cref{alg_amnet} to a series of benchmark problems in high-dimensional settings.
When presenting these results, the abbreviations ``RE'', ``SD'', ``RT'', ``vs'', and ``Iter.'' are used for ``relative error'', ``standard deviation'', ``runtime'', ``versus'', and ``Iteration'', respectively. 

Unless otherwise specified, the initial distribution for the state variable $X_0$ in \eqref{eq_sde} is uniform on two line segments $S_0$ and $S_1$ in $\R^d$, defined as
\begin{align}
    &S_0 := \bbr{s \B{1}_d: s \in [-1, 1]}, \quad \B{1}_d := (1, 1, \cdots, 1)^{\top} \in \R^{d}.\label{eq_defS2}\\
    &S_1 := \bbr{\B{l}(s): s \in [-1, 1]}, \quad \B{l}(s) := \vbr{\big}{l_1(s), l_2(s), \cdots, l_d(s)}^{\top} \in \R^d \label{eq_defS3}
\end{align}
where $l_i(s) := s\,\mathrm{sgn}\br{\sin(i)} + \cos(i + \pi s)$ and $\mathrm{sgn}(z) := -1, 0, 1$ for $z < 0$, $z = 0$, and $z > 0$, respectively.
Here $S_0$ is a line segment along the diagonal of the hypercube $[-1, 1]^d$, while $S_1$ is a more complex curve spanning the hypercube $[-2, 2]^d$.
For $i = 0, 1$, $\bar{S}_i$ denotes a uniform grid of 100 points on the line segment $S_i$.

We introduce the following indicators to quantify performance of the proposed method: 
\begin{align}
    \mathrm{RE}_{v}^{t_0} &:= \Big\{\sum_{x \in \bar{S}_0 \cup \bar{S}_1} \abs{v_{\theta}\br{0, x} - v(0, x)} \Big\} \big/ \Big\{\sum_{x \in \bar{S}_0 \cup \bar{S}_1} \abs{v(0, x)} \Big\},\\
    \mathrm{RE}_{v}^{\mathrm{path}} &:= \Big\{\sum_{n=0}^{N-1} \sum_{m=1}^{M_{\mathrm{test}}} \abs{v_{\theta}(t_n, \bar{X}_n^m) - v(t_n, \bar{X}_n^m)} \Big\} \big/ \Big\{\sum_{n=0}^{N-1} \sum_{m=1}^{M_{\mathrm{test}}} \abs{v(t_n, \bar{X}_n^m)} \Big\}, \\
    \mathrm{RE}_{\mathrm{c}} &:= \frac{\hat{J}(u_{\alpha}) - v(0, x_0)}{\abs{v(0, x_0)}}, \quad \hat{J}(u_{\alpha}) := \frac{1}{M_{\mathrm{test}}} \sum_{m=1}^{M_{\mathrm{test}}} \Big\{ \sum_{n=0}^{N-1} c_n^m \Delta t + g(\bar{X}_N^m)\Big\}. \label{eq_empJ}
\end{align}
where $c_n^m := c(t_n, \bar{X}_n^m, u_{\alpha}(t_n, \bar{X}_n^m))$, and the paths $\bar{X}_n^m$ are generated by simulating the discretized state equation~\eqref{eq_SDE} with $\hat{u} = u_{\alpha}$ and $\bar{X}_0^m = x_0 := \B{0}_d \in \R^d$ with $M_{\mathrm{test}}$ independent realizations.  
When exact solutions are unavailable, we use the reference solutions specified for the corresponding numerical examples in place of the exact value function $v$ and optimal control $u$.
We set $M_{\mathrm{test}} = 4096$ for all tests, except for $\mathrm{RE}_{v}^{\mathrm{path}}$ in subsections~\ref{sec_hjb} and~\ref{sec_complex}, for which we use $M_{\mathrm{test}} = 128$ because computing the reference solutions on the paths $\bar{X}_n^m$ is prohibitively expensive.
The indicator $\mathrm{RE}_{u}^{\mathrm{path}}$ is defined analogously to $\mathrm{RE}_{v}^{\mathrm{path}}$, with $v_{\theta}$ replaced by $u_{\alpha}$ and $v$ by $u$.

The above indicators can be interpreted as follows. 
The indicator $\mathrm{RE}_{v}^{t_0}$ quantifies the accuracy of the learned value function $v_{\theta}(t, x)$ at $t = 0$ for $x$ on $S_0$ and $S_1$.
The line segments $S_0$ and $S_1$ are used to better visualize the solution profile in high dimensions. 
The indicators $\mathrm{RE}_{v}^{\mathrm{path}}$ and $\mathrm{RE}_{u}^{\mathrm{path}}$ quantify the accuracy of $v_{\theta}$ and $u_{\alpha}$, respectively, along the simulated paths $\bar{X}_n^m$.
They thus assess the accuracy of our method in regions of the high-dimensional state space explored by the dynamics.
The indicator $\mathrm{RE}_{\mathrm{c}}$ denotes the relative error of the cost functional $J(u_{\alpha})$ in \eqref{eq_costfunc}, evaluated at the learned control $u_{\alpha}$. 
Ideally, if $u_{\alpha}$ is the optimal control and there is no Monte Carlo error, $\mathrm{RE}_{\mathrm{c}}$ is zero; thus, it measures the optimality of the learned control $u_{\alpha}$.

\subsection{Parameter settings}\label{sec_param}

The following parameter settings serve as the default choices for \Cref{alg_amnet} in the numerical experiments.
For specific examples, certain parameters are adjusted either to ensure fair comparisons with existing methods or to test the scalability of the proposed method; the corresponding values are reported in the respective subsections.

The neural networks $u_{\alpha}$ and $v_{\theta}$ are constructed as in \eqref{eq_utx} and \eqref{eq_vtx}, respectively. 
Their inner networks $\psi_{\alpha}$ and $\phi_{\theta}$ are fully connected with $H = 4$ hidden layers, each containing $W$ ReLU units; the values of $W$ used in the experiments are reported together with the corresponding numerical results. 
The adversarial network $\rho_{\eta}$ is taken as in \eqref{eq_defrhonet}, with output dimension $r = 600$ and scaling parameter $c = 100$, except for the example in \Cref{sec_complex}, which uses $c = 10$.

For parameter settings of \Cref{alg_amnet}, the total number $M$ of sample paths retained in memory is set to $M = 10^4$ for all the tests, except for the scalar-channel example in \Cref{sec_k1pair}, which uses $M = 10{,}240$.
The number of time steps in \eqref{eq_timepart} is set to $N = 32$ unless otherwise specified; the comparison with SOC-MartNet \cite{cai2024socmartnet} in \Cref{sec_hjb} adopts $N = 100$ to align the time grid with that of SOC-MartNet.
The mini-batch sizes $\abs{\mathbb{A}_1} = \abs{\mathbb{A}_2}$ are set to $200$, except for the example in \Cref{sec_complex}, which uses $\abs{\mathbb{A}_1} = \abs{\mathbb{A}_2} = 400$.
The interval for updating pilot paths and the percentage of paths updated are set to $I_0 = M/\abs{\mathbb{A}_1 \cup \mathbb{A}_2}$ and $r\% = 20\%$, respectively.
The inner iteration steps are $J = 2K = 2$.
The stochastic gradient updates are performed by the RMSProp algorithm. 
At the $i$-th iteration step, the learning rates of $u_{\alpha}$ and $v_{\theta}$ are both set to $\delta_0 \times 10^{-3} \times 0.01^{i/I}$, and the learning rate of $\rho_{\eta}$ is set to $\delta_0 \times 10^{-2} \times 0.01^{i/I}$, where $\delta_0 := 3d^{-0.5}$ for $d \leq 1000$ and $\delta_0 := 3 d^{-0.8}$ for $d > 1000$.

All experiments are implemented in Python 3.12 using PyTorch 2.5.1, where bfloat16 automatic mixed precision is used for the forward and backward passes of the neural networks and float32 precision for all error evaluations.
The algorithm is accelerated by the strategy of Distributed Data Parallel (DDP) \footnote{\url{https://github.com/pytorch/tutorials/blob/main/intermediate_source/ddp_tutorial.rst}} on compute nodes equipped with NVIDIA A100-SXM4-80GB GPUs, with the number of GPUs per run reported together with the corresponding results.

\subsection{HJB equations}\label{sec_hjb}

\begin{table}[t]
    \centering
    \caption{
    Results of \Cref{alg_amnet} for HJB-1, -2, and -3 in \Cref{sec_hjb}. Mean and SD denote the mean and standard deviation of each reported value over five independent runs.
    }\label{tab_hjb1e4}
    \footnotesize
    \setlength{\tabcolsep}{2pt}
    \resizebox{\textwidth}{!}{%
    \begin{tabular}{c c c c c c c c c c c}
    \toprule
    Equation & $d$ & $W$ & $H$ & \multicolumn{2}{c}{$\mathrm{RE}_{v}^{t_0}$} & \multicolumn{2}{c}{$\mathrm{RE}_{v}^{\mathrm{path}}$} & \multicolumn{2}{c}{$\mathrm{RE}_{\mathrm{c}}$} & RT (s) \\
    \cmidrule(lr){5-6} \cmidrule(lr){7-8} \cmidrule(lr){9-10}
     & & & & Mean & SD & Mean & SD & Mean & SD &  \\
    \midrule
    HJB-1 & $10{,}000$ & $d+10$ & 4 &  2.44E-3 & 1.8E-4 & 9.2E-3 & 7.9E-4 & 1.3E-3 & 1.7E-4 & 2812 \\
    HJB-2 & $10{,}000$ & $d+10$ & 4 & 2.03E-2 & 2.0E-3 & 6.2E-3 & 2.1E-3 & 8.3E-3 & 2.4E-4 & 2814 \\
    HJB-3 & $10{,}000$ & $d+10$ & 6 & 3.69E-2 & 2.6E-3 & 2.1E-2 & 2.5E-3 & 1.3E-3 & 5.2E-5 & 3953 \\
    HJB-3 & $2{,}000$ & $d+10$ & 4 & 6.24E-2 & 6.4E-3 & 4.1E-2 & 7.9E-3 & 2.8E-2 & 1.0E-3 & 248 \\
    HJB-3 & $2{,}000$ & $5d+10$ & 4 & 3.20E-2 & 4.7E-3 & 3.5E-2 & 2.4E-3 & 7.6E-4 & 6.5E-5 & 1791 \\
    \bottomrule
    \end{tabular}}
\end{table}

\begin{figure}[t]
    \centering
    \subfloat[Solution profile, HJB-1]{\includegraphics[width=0.3\textwidth]{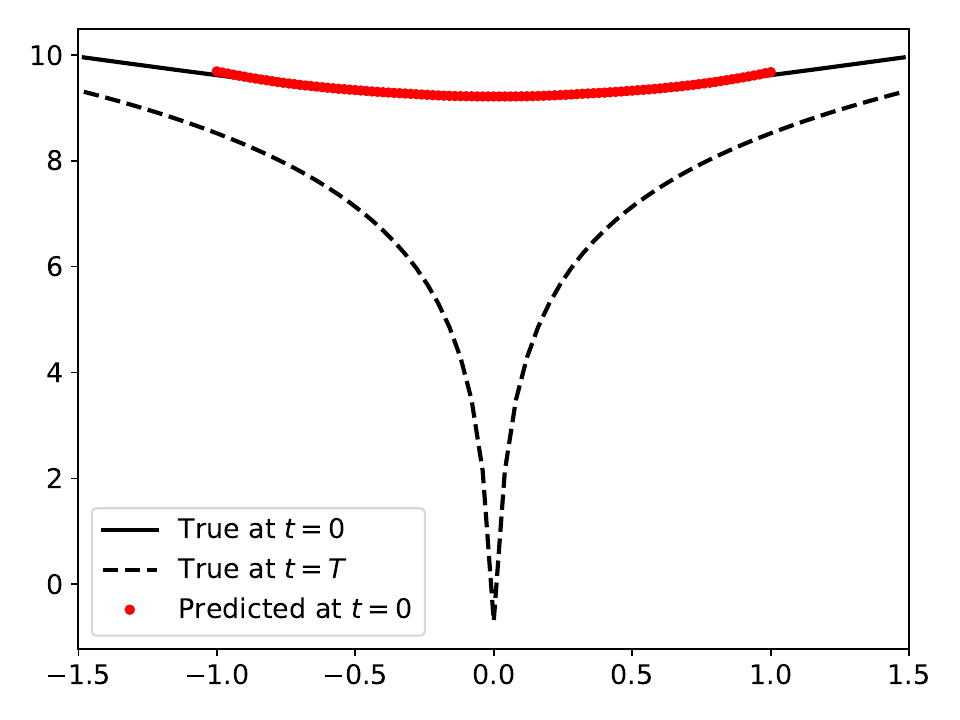}}
    \subfloat[Solution profile, HJB-2]{\includegraphics[width=0.3\textwidth]{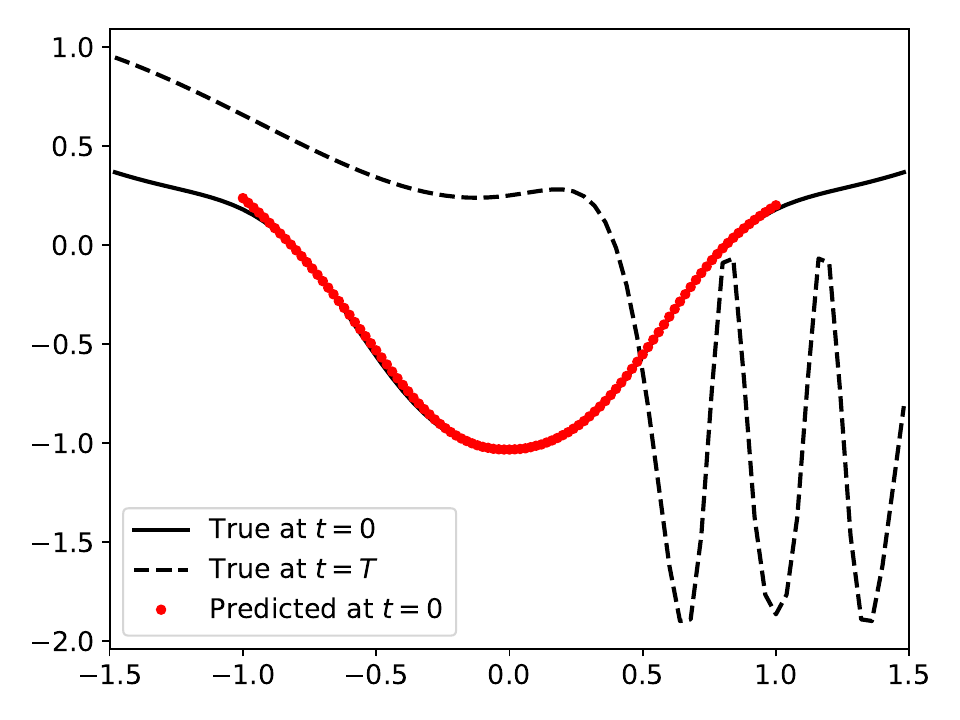}}
    \subfloat[Solution profile, HJB-3]{\includegraphics[width=0.3\textwidth]{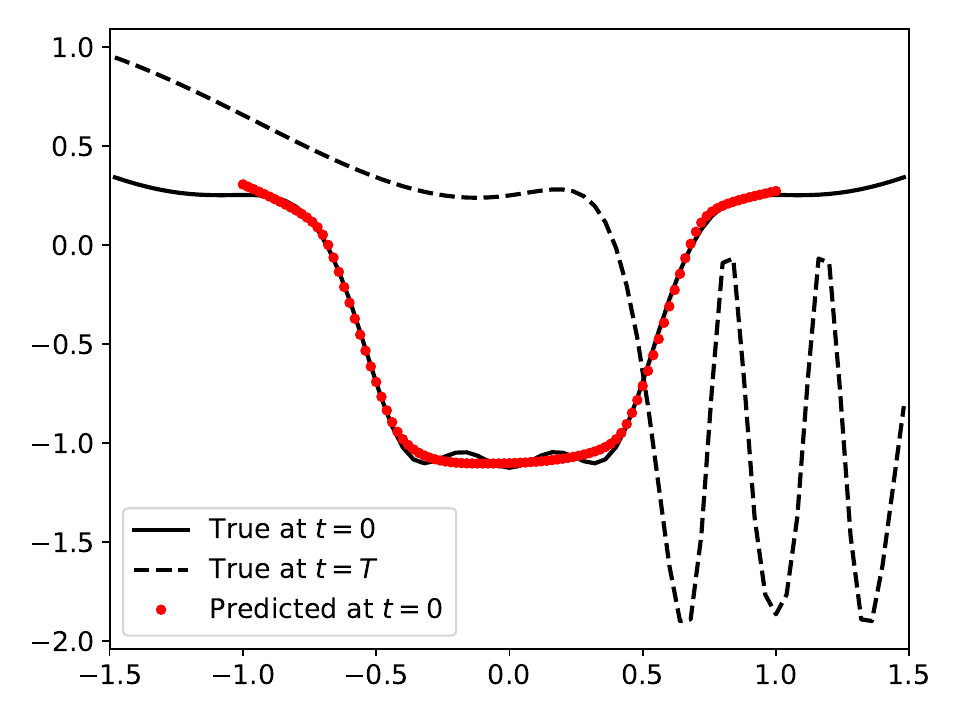}}
    \\
    \subfloat[RE vs Iter., HJB-1]{\includegraphics[width=0.3\textwidth]{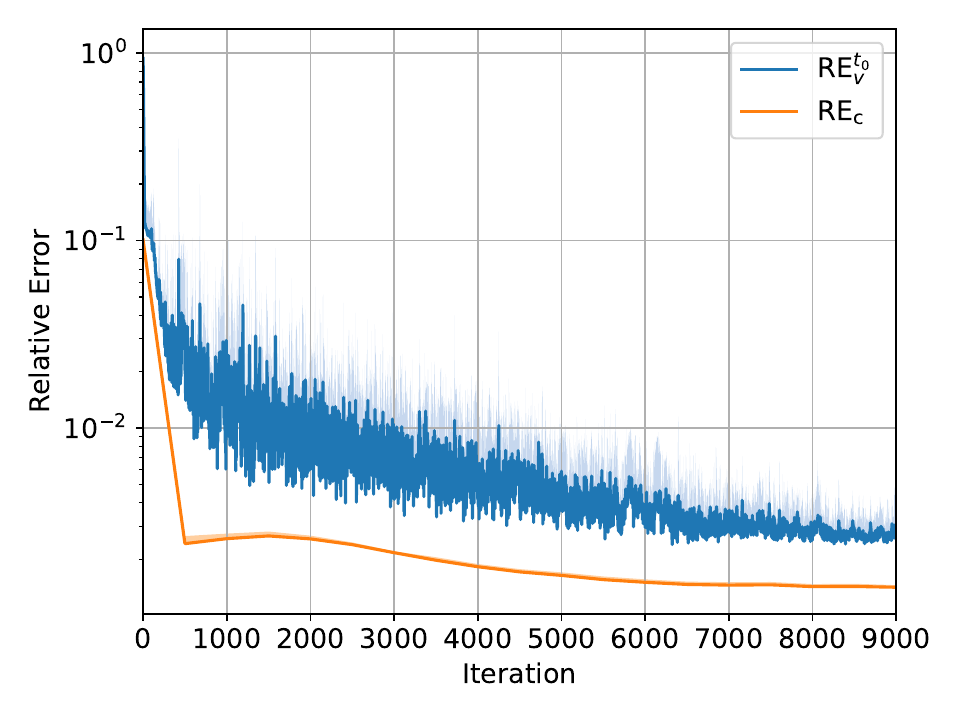}}
    \subfloat[RE vs Iter., HJB-2]{\includegraphics[width=0.3\textwidth]{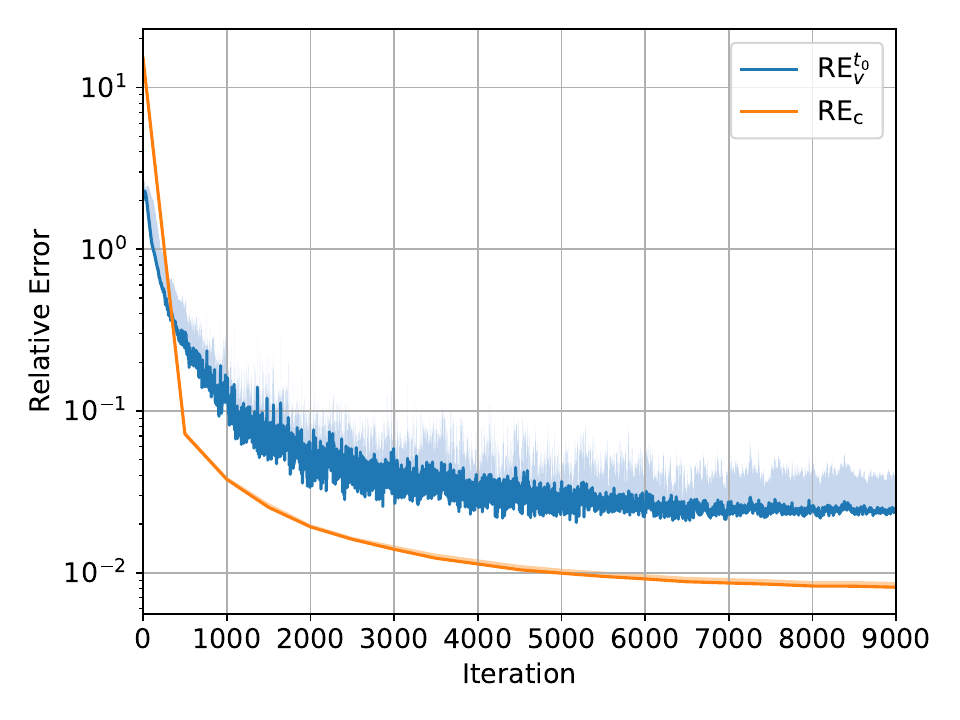}}
    \subfloat[RE vs Iter., HJB-3]{\includegraphics[width=0.3\textwidth]{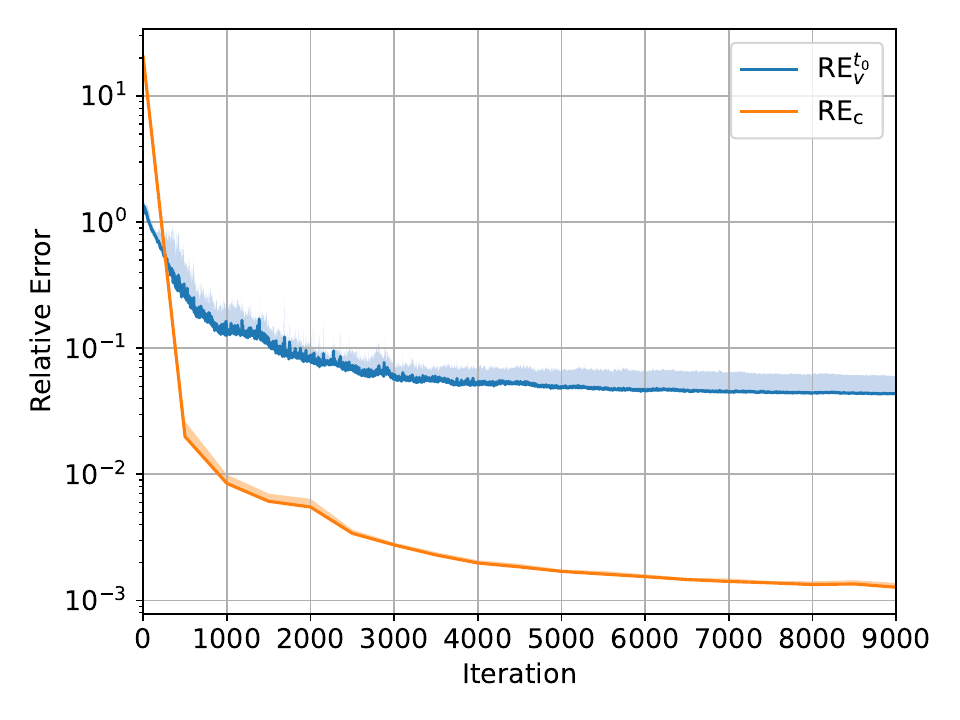}}
    \caption{
    Results for HJB-1, -2, and -3 in \Cref{sec_hjb} with $d = 10^4$. Shading denotes the five-run mean plus or minus twice the SD in (a)--(c) and the five-run mean plus twice the SD in (d)--(f). The final-iteration statistics are reported in \Cref{tab_hjb1e4}.
    }\label{fig1_d1e4}
\end{figure}

\begin{figure}[t]
    \centering
    \subfloat[$s \mapsto v_{\theta}(0, s\B{1}_d) - v(0, s\B{1}_d)$]{\includegraphics[width=0.32\textwidth]{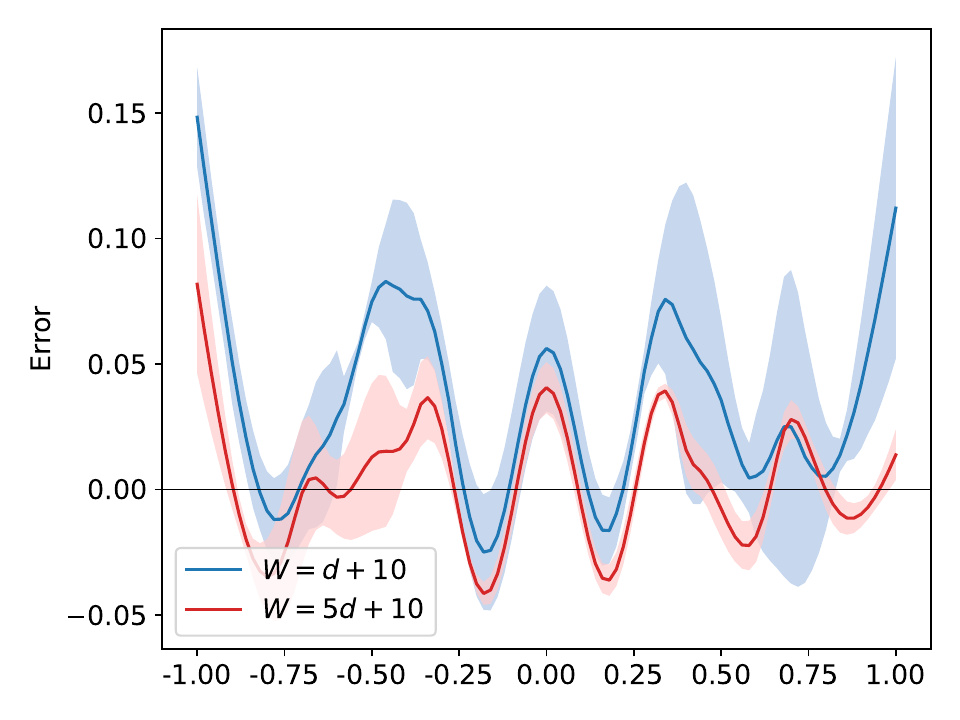}}
    \subfloat[$\mathrm{RE}_{v}^{t_0}$ vs Iter.]{\includegraphics[width=0.32\textwidth]{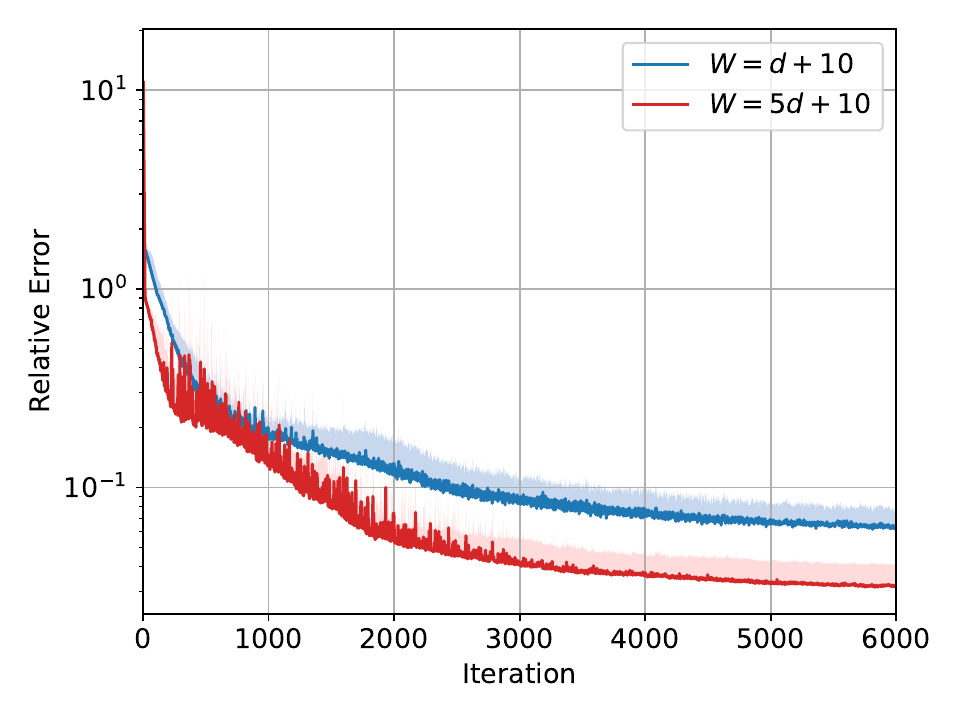}}
    \subfloat[$\mathrm{RE}_{\mathrm{c}}$ vs Iter.]{\includegraphics[width=0.32\textwidth]{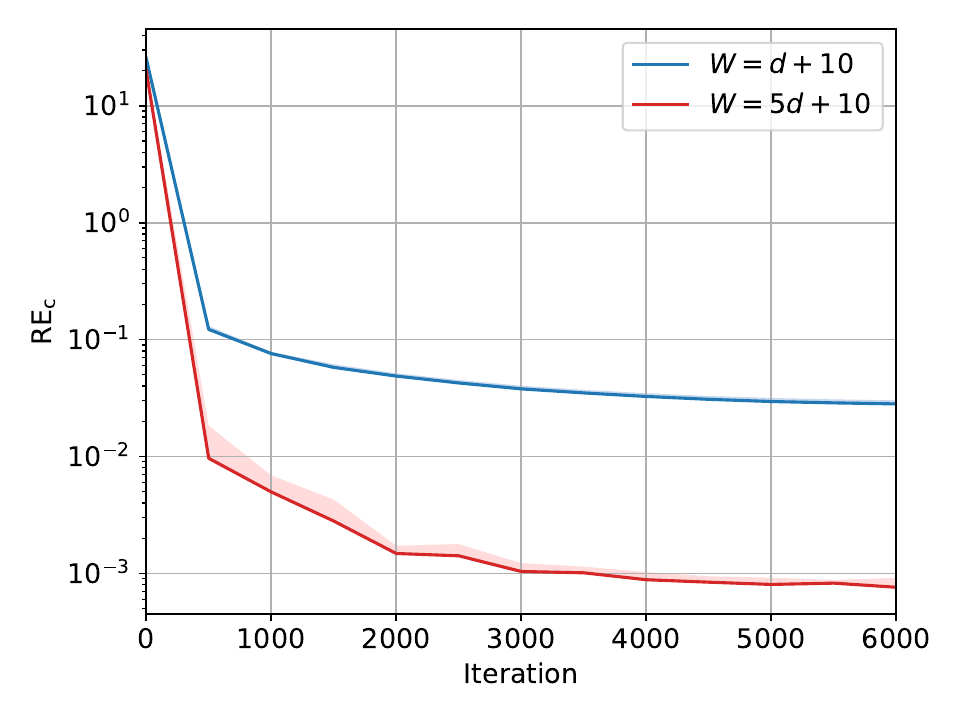}}
    \caption{Results for HJB-3 in \Cref{sec_hjb} with $d = 2000$ and $W = d + 10$ or $W = 5d + 10$. Shading denotes the five-run mean plus or minus twice the SD in (a) and the five-run mean plus twice the SD in (b)--(c). The final-iteration statistics are reported in \Cref{tab_hjb1e4}.
    }\label{fig_HJB3b}
\end{figure}

\begin{figure}[t]
    \centering
    \subfloat[HJB-2]{\includegraphics[width=0.41\textwidth]{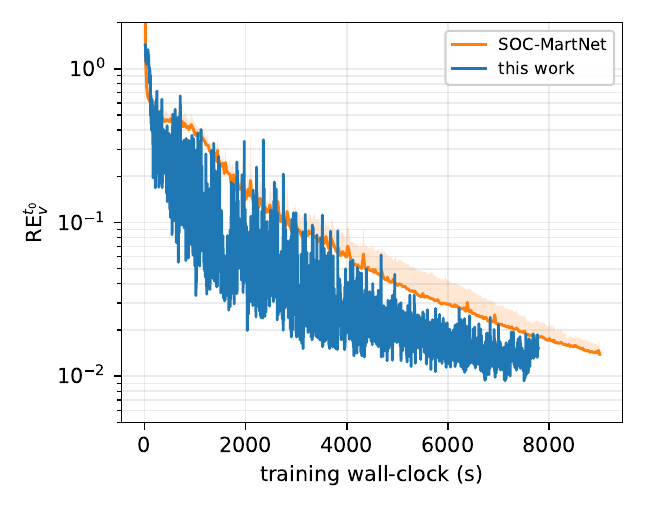}}
    \subfloat[HJB-3]{\includegraphics[width=0.41\textwidth]{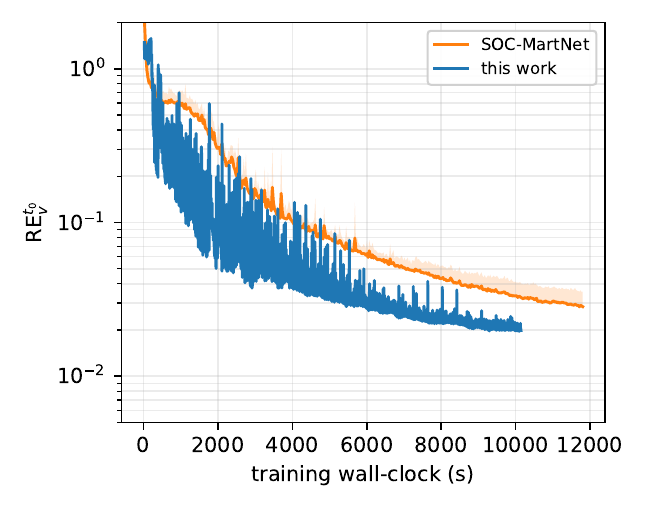}}
    \caption{Relative errors versus training wall-clock time for \Cref{alg_amnet} and SOC-MartNet \cite{cai2024socmartnet} on HJB-2 and -3 in \Cref{sec_hjb} with $d = 10^4$. The SOC-MartNet curves and shading denote the five-run mean and the mean plus twice the SD, respectively.}\label{fig_cmp_socmartnet}
\end{figure}

The considered HJB equations are given by \eqref{eq_HJBPDE} with
\begin{equation}\label{eq2_hjbtest}
    \mathcal{L}^{\kappa} = (b + 2 \kappa)^{\top} \partial_x + \delta^2 \mathrm{Tr}[\partial_{xx}^2], \quad c(t, x, \kappa) = \delta^{-2} \abs{\kappa}^2, \quad U = \R^d, \quad T = 1,
\end{equation}
where $\delta > 0$, $b \in \R^d$ and the terminal condition $v(T, x) = g(x)$ will be specified later. 
By the It\^o formula, the resulting HJB \eqref{eq_HJBPDE} admits an analytic solution given by  
\begin{equation}\label{eq_vtx_hjb1}
    v(t, x) = - \ln\br{\E{\exp \br{-g(x + (T-t)b + \sqrt{2} \delta B_{T-t})}}}
\end{equation}
for $(t, x) \in [0, T) \times \R^d$.
For a reference solution, the expectation in \eqref{eq_vtx_hjb1} is approximated by the Monte Carlo method using $10^6$ i.i.d.\ samples of $B_{T-t}$.

We consider the following specific HJB equations, which are given by \eqref{eq_HJBPDE} and \eqref{eq2_hjbtest} with different parameters: 
\begin{itemize}
    \item HJB-1: $b = 0$, $\delta = 1$, $v(T, x) = \ln(0.5 (1+\abs{x}^2))$. 
    
    \item HJB-2: $b = \B{1}_d$, $\delta = 0.2$, and $v(T, x) = \bar{g}(x - \B{1}_d)$ with
    \begin{equation*}
        \bar{g}(x) := \frac{1}{d} \sum_{i=1}^d \bbr{\sin(x_i - \frac{\pi}{2}) + \sin\fbr{\fbr{0.3/\pi + x_i^2}^{-1}}}. 
    \end{equation*}
    \item HJB-3: a variant of HJB-2 with $\delta = 0.1$. 
\end{itemize}

We have the following comments for the above equations:
HJB-1 is a benchmark problem commonly used in existing works, such as \cite{weinan2017deep,Ji2022Solving,Bachouch2022Deep,han2018solving,wang2022is,He2023Learning,raissi2018forwardbackward,hu2024sdgd}.
HJB-2 and HJB-3 are more challenging than HJB-1 due to their terminal function $\bar{g}(x)$, which is highly oscillatory around $x = 0$, resulting in an uneven solution when $\delta$ is small; see \Cref{fig1_d1e4}(b) and (c).

By applying \Cref{alg_amnet} to the above HJB equations, the relevant numerical results are presented in \Cref{fig_HJB3b,fig1_d1e4} and \Cref{tab_hjb1e4}.
The results in \Cref{fig_HJB3b} and \Cref{tab_hjb1e4} provide evidence that \Cref{alg_amnet} remains effective for HJB equations with state dimension up to $10^4$, including the oscillatory HJB-2 and HJB-3 cases. The small value-function errors on both the prescribed slices and simulated paths, together with the small cost errors, indicate that the method learns not only an accurate value function but also a feedback control with near-optimal performance on the explored state-space region. 
As shown in \Cref{fig_HJB3b}, the improvement obtained by widening the network for HJB-3 suggests that, for problems with more complicated solution profiles, network expressivity rather than the martingale formulation itself can become the practical accuracy bottleneck.

A comparison with SOC-MartNet \cite{cai2024socmartnet} on HJB-2 and HJB-3 is presented in \Cref{fig_cmp_socmartnet}.
For both HJB-2 and HJB-3, \Cref{fig_cmp_socmartnet} shows that the proposed method attains lower values of $\mathrm{RE}_{v}^{t_0}$ for most of the common range of training wall-clock times.
This advantage is consistent with the derivative-free formulation: SOC-MartNet relies on automatic differentiation to evaluate the first-order spatial derivatives of $v_{\theta}$, whereas the proposed method completely avoids derivative evaluations.

\subsection{SOCP with complex state dynamics}\label{sec_complex}

% \begin{figure}[t]
%     \centering
%     \subfloat[$v(0, x)$ for $x \in S_0 \cup S_1$]{\includegraphics[width=0.4\textwidth]{Figures/cplxd1e4_profiles_S0S1.pdf}}
%     \subfloat[$\mathrm{RE}_{v}^{t_0}$ and $\mathrm{RE}_{\mathrm{c}}$ vs Iter.]{\includegraphics[width=0.4\textwidth]{Figures/cplxd1e4_rerc_hist.pdf}}
%     \caption{Results for the complex-dynamics SOCP in \Cref{sec_complex} with $d = 10^4$. Shading denotes the five-run mean plus or minus twice the SD in (a) and the five-run mean plus twice the SD in (b).}\label{fig_cplx}
% \end{figure}

% \begin{figure}[t]
%     \centering
%     \begin{resizebox}{0.8\textwidth}{!}
%         {\input{fig_scaling_tikz}}
%     \end{resizebox}
%     \caption{Sensitivity to the test-function output dimension $r$ for the SOCP example in \Cref{sec_complex} with $d = 10^4$. Shading denotes the five-run mean plus twice the SD.}\label{fig_scaling}
% \end{figure}

% \begin{figure}[t]
%     \centering
%     \input{fig_cost_tikz}
%     \caption{Relative per-iteration runtime and GPU memory under one-at-a-time parameter sweeps around the configuration of \Cref{sec_complex}. Values are normalized by those at the smallest displayed value of the varying parameter. Each point uses one four-GPU run, and runtime is the median over iterations 51--200.}\label{fig_cost}
% \end{figure}

Now we consider a more complex SOCP, in which the components of the state variable $X$ are asymmetrically coupled, to assess the performance of \Cref{alg_amnet} for problems with intricate dynamics.
The problem is an instance of the SOCP~\eqref{eq_SOCP} with the drift and diffusion coefficients in the state equation~\eqref{eq_sde} specified as $\mu(t, x, \kappa) = b + c_0 s \kappa$ and $\sigma(t, x, \kappa) = s$ for $(t, x, \kappa) \in [0, T] \times \R^d \times \R^d$, where $b$ and $s$ are functions of $(t, x)$, defined by
\begin{equation*}
    b_i = 2 \sin \bigl(t + i + \sum_{k=1}^d 0.9^k x_{i+k} - 1 \bigr), \quad s_{ij} = \sqrt{2} c_1 0.9^{\abs{i-j}} \bigl(1 + 0.5 \cos(t + i + x_{i+1})\bigr)
\end{equation*}
for $i, j = 1, 2, \cdots, d$ with $x_{d+k} := x_k$, $c_0 = 2$ and $c_1 = 0.5$.
The control set is given by $U = \R^d$.
The cost functional is given by \eqref{eq_costfunc} with 
\begin{equation}
    c(t, x, \kappa) = \frac{1}{2} \abs{\kappa}^2, \quad g(x) = \ln\vbr{\Big}{1 + \frac{1}{d} \sum_{i=1}^d x_i^2} + \frac{0.3}{d} \sum_{i=1}^d \sin\br{0.5 x_i} + \frac{2}{d} \sum_{i=1}^d \sin x_i. 
\end{equation}
By the Cole-Hopf transformation, the value function admits the representation \cite[Appendix C]{Cai2026Deep}
\begin{equation}\label{eq_vtx_cplx}
    v(t, x) = - 0.25 \ln \E{\exp\br{-4 g\br{X_T^{0}}} \big\vert X_t^{0} = x}, \quad t \in [0, T), \quad x \in \R^d,
\end{equation}
where $X^{0}$ refers to the state process $X$ given by \eqref{eq_sde} with $u \equiv 0$. 
The reference solution is computed by approximating the expectation in \eqref{eq_vtx_cplx} via Monte Carlo simulation with $10^6$ samples of $X_T^{0}$, generated by the Euler-Maruyama discretization of \eqref{eq_sde} with $N = 100$ time steps.

\begin{figure}[t]
    \centering
    \subfloat[$v(0, x)$ for $x \in S_0 \cup S_1$]{\includegraphics[width=0.4\textwidth]{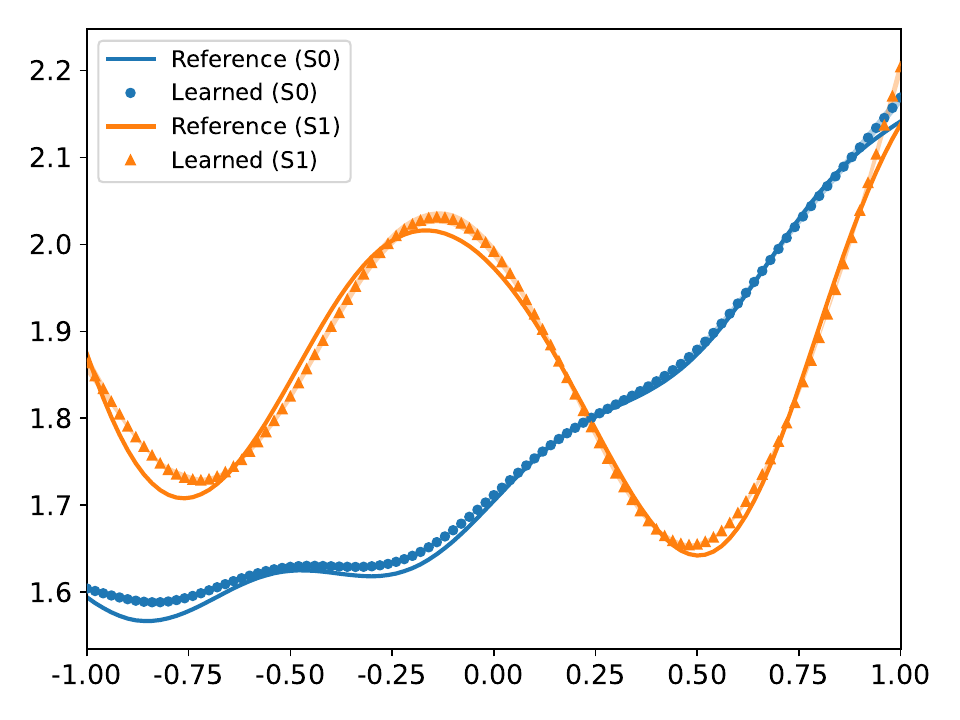}}
    \subfloat[$\mathrm{RE}_{v}^{t_0}$ and $\mathrm{RE}_{\mathrm{c}}$ vs Iter.]{\includegraphics[width=0.4\textwidth]{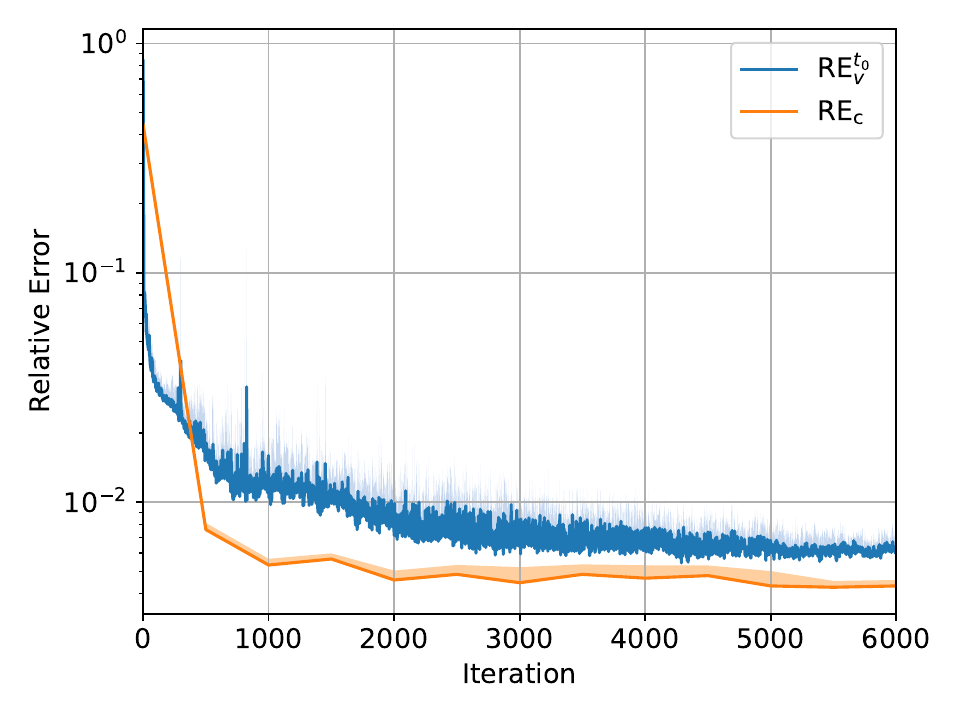}}
    \caption{Results for the complex-dynamics SOCP in \Cref{sec_complex} with $d = 10^4$. Shading denotes the five-run mean plus or minus twice the SD in (a) and the five-run mean plus twice the SD in (b).}\label{fig_cplx}
\end{figure}

\Cref{fig_cplx} presents the results of \Cref{alg_amnet} over 5 independent runs.
Despite the asymmetric coupling and nonlinear state dependence, the learned profiles closely track the reference solution on both $S_0$ and $S_1$. 
The final values of $\mathrm{RE}_{v}^{t_0}$ and $\mathrm{RE}_{\mathrm{c}}$, both below $10^{-2}$, indicate that the learned value function is accurate and that the learned feedback control achieves near-optimal performance, respectively.

In the following, we examine the sensitivity of \Cref{alg_amnet} to the output dimension $r$ of the test function network $\rho$ for this SOCP example.
The results are reported in \Cref{fig_scaling}.
As shown in \Cref{fig_scaling}, increasing $r$ from $75$ to $1200$ reduces $\mathrm{RE}_{v}^{t_0}$ by approximately $1/3$, with diminishing marginal improvement as $r$ increases. 
Meanwhile, the runtime and GPU-memory ratios vary by less than $1.5\%$ over the tested range. 
Thus, \Cref{alg_amnet} is robust with respect to the choice of $r$: increasing $r$ offers a modest accuracy improvement without a significant additional computational cost.

Figure~\ref{fig_cost} reports the per-iteration runtime and GPU memory as the time-step number $N$, total number $M$ of pilot paths, state dimension $d$, and mini-batch size $\abs{\mathbb{A}_1} = \abs{\mathbb{A}_2}$ vary one at a time, with all other parameters fixed at their default values. Increasing $N$ from $2^3$ to $2^7$ changes the runtime by less than $6\%$ but increases GPU memory by a factor of approximately $2.6$, consistent with temporal parallelization. Increasing the mini-batch size mainly raises the runtime, whereas increasing the pilot-pool size $M$ mainly raises the memory required to retain pilot paths. 
The state dimension $d$ is the dominant bottleneck: when $d$ increases from $2^{10}$ to $2^{14}$, runtime and GPU memory increase by factors of approximately $34$ and $44$, respectively. 
Both factors are well below the factor of $2^8 = 256$ implied by quadratic scaling;
thus, over the tested range, the observed resource growth is subquadratic rather than exponential and provides no empirical indication of a curse of dimensionality.

\begin{figure}[t]
    \centering
    \begin{resizebox}{0.8\textwidth}{!}
        {\input{fig_scaling_tikz}}
    \end{resizebox}
    \caption{Sensitivity to the test-function output dimension $r$ for the SOCP example in \Cref{sec_complex} with $d = 10^4$. Shading denotes the five-run mean plus twice the SD.}\label{fig_scaling}
\end{figure}

\begin{figure}[t]
    \centering
    \input{fig_cost_tikz}
    \caption{Relative per-iteration runtime and GPU memory under one-at-a-time parameter sweeps around the configuration of \Cref{sec_complex}. Values are normalized by those at the smallest displayed value of the varying parameter. Each point uses one four-GPU run, and runtime is the median over iterations 51--200.}\label{fig_cost}
\end{figure}

\subsection{SOCP with a scalar control channel}\label{sec_k1pair}

% \begin{figure}[t]
%     \centering
%     \subfloat[REs vs Iter.]{\includegraphics[width=0.32\textwidth]{Figures/fig_k1x0pt_hist.pdf}}
%     \subfloat[QQ plot of $u_{\alpha}$ vs $u^{*}$]{\includegraphics[width=0.32\textwidth]{Figures/k1x0_uqq_paths.pdf}}
%     \subfloat[QQ plot of $v_{\theta}$ vs $v$]{\includegraphics[width=0.32\textwidth]{Figures/k1x0_vqq_paths.pdf}}
%     \caption{Results for the SOCP in \Cref{sec_k1pair} with $d = 10^4$. 
%     In (a), shading denotes the five-run mean plus twice the SD.
%     In (b) and (c), the quantile-quantile (QQ) plots compare $u_{\alpha}$ and $v_{\theta}$, respectively, evaluated on the points $(t_n, \bar{X}_n^m)$ introduced in \eqref{eq_empJ}.
%     }\label{fig_k1pair}
% \end{figure}

We consider a SOCP with a scalar control channel to show the effectiveness of \Cref{alg_amnet} in approximating the optimal control in underactuated problems.
The problem is given by \eqref{eq_SOCP} with the drift and diffusion coefficients in the state equation~\eqref{eq_sde} specified as
\begin{align*}
    &\mu_i(t, x, \kappa) = x_{i+1} - x_i + 2 q_i \kappa, \quad \sigma(t, x, \kappa) = \sqrt{2} Q Q^{\top}, \quad (t, x, \kappa) \in [0, T] \times \R^d \times U,\\
    &x_{d+1} := x_1, \quad q := (q_i)_{i=1}^d = 0.5 \times (\sqrt{2}, \, \sqrt{2}, \,0, \cdots, 0)^{\top} \in \R^d, \quad T = 1 
\end{align*}
for $i = 1, 2, \cdots, d$, where $X_0 = \B{0}_d$ and $U = \R$.
Here the columns of $Q \in \R^{d \times 20}$ are the first $20$ vectors obtained by Gram--Schmidt orthogonalization of the sequence $q, \varphi^{(1)}, \cdots, \varphi^{(d)}$, where $\varphi^{(1)}, \cdots, \varphi^{(d)}$ is the cosine basis of $\R^d$, with the $j$-th entry of the $k$-th vector given by
\begin{equation*}
    \varphi^{(k)}_j = a_k \cos\Bigl( (2d)^{-1} \pi (2 j - 1)(k - 1) \Bigr), \quad j, k = 1, 2, \cdots, d,
\end{equation*}
with $a_1 := 1/\sqrt{d}$ and $a_k := \sqrt{2 / d}$ for $k \ge 2$.
The cost functional is given by \eqref{eq_costfunc} with
\begin{equation}\label{eq_cost_k1}
    c(t, x, \kappa) = \kappa^2, \quad g(x) = 0.5 \br{q^{\top} x - 1}^2. 
\end{equation}
The optimal feedback control and the value function are given by
\begin{equation}\label{eq_kstar_k1}
    u^{*}(t, x) = - q^{\top} \nabla v(t, x), \quad v(t, x) = \frac{1}{2} \log\big(1 + C_{T-t}\big) + \frac{\br{q^{\top} e^{(T-t) L} x - 1}^2}{2 \br{1 + C_{T-t}}},
\end{equation}
where $C_{T-t} := 2 \int_0^{T-t} \abs{q^{\top} e^{r L} Q}^2 \di r$ and the matrix $L \in \R^{d \times d}$ is such that, for every $x \in \R^d$, the $i$-th element of $L x$ equals $x_{i+1} - x_i$ for $i = 1, 2, \cdots, d$, with $x_{d+1} := x_1$.

\begin{figure}[t]
    \centering
    \subfloat[REs vs Iter.]{\includegraphics[width=0.32\textwidth]{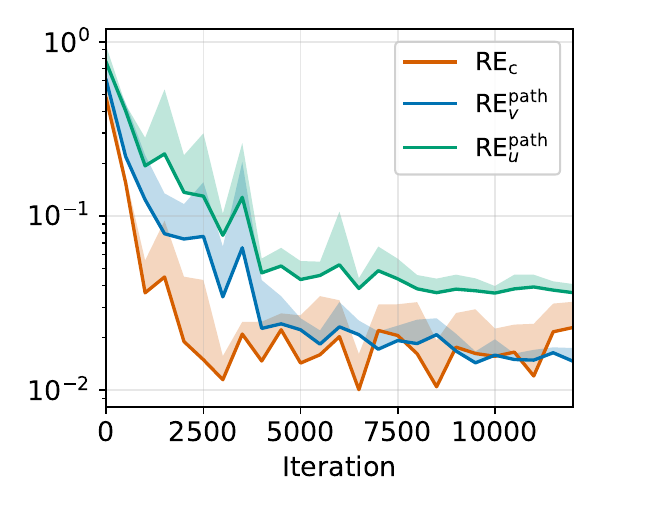}}
    \subfloat[QQ plot of $u_{\alpha}$ vs $u^{*}$]{\includegraphics[width=0.32\textwidth]{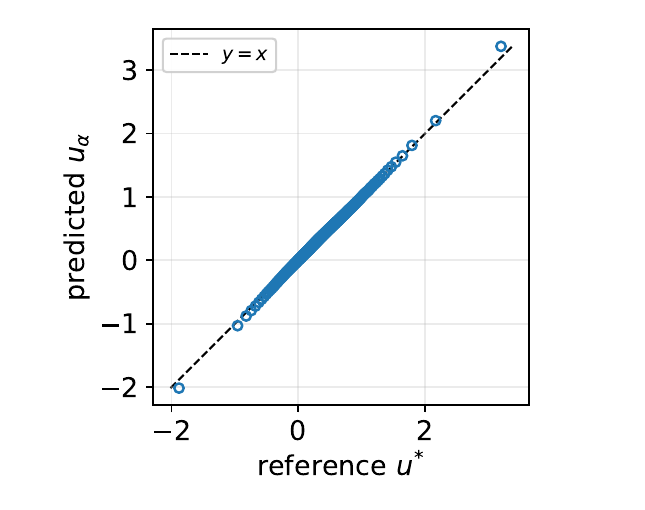}}
    \subfloat[QQ plot of $v_{\theta}$ vs $v$]{\includegraphics[width=0.32\textwidth]{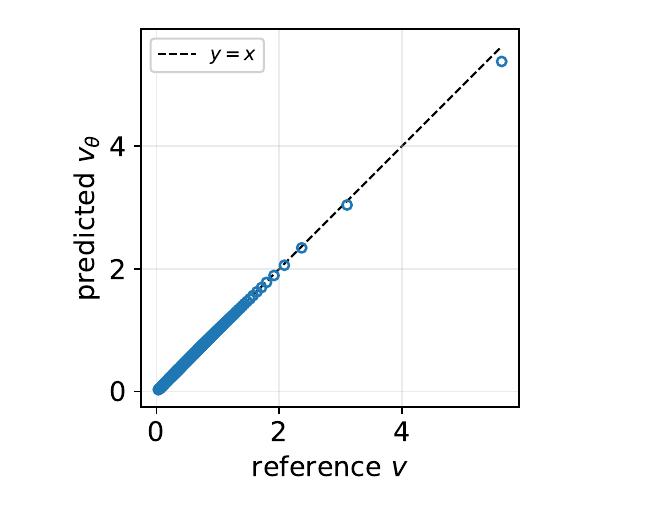}}
    \caption{Results for the SOCP in \Cref{sec_k1pair} with $d = 10^4$. 
    In (a), shading denotes the five-run mean plus twice the SD.
    In (b) and (c), the quantile-quantile (QQ) plots compare $u_{\alpha}$ and $v_{\theta}$, respectively, evaluated on the points $(t_n, \bar{X}_n^m)$ introduced in \eqref{eq_empJ}.
    }\label{fig_k1pair}
\end{figure}

The results for $d = 10^4$ are shown in \Cref{fig_k1pair}. 
The final values of $\mathrm{RE}_{\mathrm{c}}$, $\mathrm{RE}_{v}^{\mathrm{path}}$, and $\mathrm{RE}_{u}^{\mathrm{path}}$ are all below $4 \%$. 
In contrast to evaluating the value function only on the spatial curves $S_0$ and $S_1$ at $t = 0$, this example evaluates both $v_{\theta}$ and $u_{\alpha}$ along the paths driven by the learned control. 
The close agreement with the diagonal in the two QQ plots further confirms that both approximations remain accurate throughout this path-explored region in high-dimensional state space.

\subsection{HJB equation without explicit \texorpdfstring{$\boldsymbol{\inf_u H}$}{inf\_u H}}\label{sec_eps}

Now we consider a SOCP whose complexity arises in the Hamiltonian, for which the associated HJB equation reads:
\begin{equation}\label{eq_hjb_eps}
    \partial_t v + \inf_{\kappa \in \R^d} \bbr{(b + \kappa)^{\top} \partial_x v + \frac{1}{2}\abs{\kappa}^2 + \epsilon \sin(\B{1}_d^{\top} \kappa)} + \frac{1}{8} \Delta_x v = 0 
\end{equation}
for $(t, x) \in [0, T) \times \R^d$, where $\epsilon$ is a scalar parameter specified later, and $b$ is a $\R^d$-valued function of $(t, x)$ given by
\begin{equation*}
    b_i = \sin(t + i + x_{i+1} - 1), \quad x_{d+1} := x_1, \quad i = 1, 2, \cdots, d,
\end{equation*}
The terminal condition is
\begin{equation*}
    v(T, x) = g(x) := \ln\Big(1 + \frac{1}{d} \sum_{i=1}^d \br{x_i^2 + 0.5 \sin(10 x_i)}\Big),\quad x \in \R^d, \quad T = 1.
\end{equation*}
This HJB equation arises from the following SOCP: 
\begin{equation*}
    \min_{u \in \mathcal{U}_{\mathrm{ad}}} J(u), \quad J(u) := \int_0^T \frac{1}{2}\E{\abs{u(t, X_t^u)}^2 + \epsilon \sin(\B{1}_d^{\top} u(t, X_t^u))} \di t + \E{g(X_T^u)},
\end{equation*}
where $X^u: [0, T] \times \Omega \to \R^d$ is the state process defined by 
\begin{equation}\label{eq_xtu}
    X_t^u = X_0^u + \int_0^t \br{b + u}(s, X_s^u) \di s + \int_0^t \frac{1}{2} \di B_s, \quad t \in [0, T]. 
\end{equation}
The admissible control set $\mathcal{U}_{\mathrm{ad}}$ consists of all Borel measurable functions from $[0, T] \times \R^d$ to $\R^d$, such that the SDE~\eqref{eq_xtu} admits a unique strong solution.

When $\epsilon = 0$, the value function $v(t, x)$ admits a closed-form solution given by \eqref{eq_vtx_cplx}, where $X_t^{0}$ is replaced by the solution of \eqref{eq_xtu} with $u \equiv 0$.
To compute the reference solution, we approximate the expectation in \eqref{eq_vtx_cplx} by Monte Carlo, generating $10^6$ samples of $X_T^{0}$ via the Euler--Maruyama scheme applied to \eqref{eq_xtu}.
For $\epsilon \neq 0$, the infimum in the HJB equation \cref{eq_hjb_eps} does not admit a closed-form solution due to the presence of the sine term.
Therefore, this example serves as a test case for evaluating the effectiveness of our method on general problems.
In our experiments, we consider a decreasing sequence of the perturbation parameter $\epsilon = 1, 1/2, 1/4, 1/8, 0$ with $d = 10^4$.

\begin{figure}[t]
    \centering
    \includegraphics[width=0.54\textwidth]{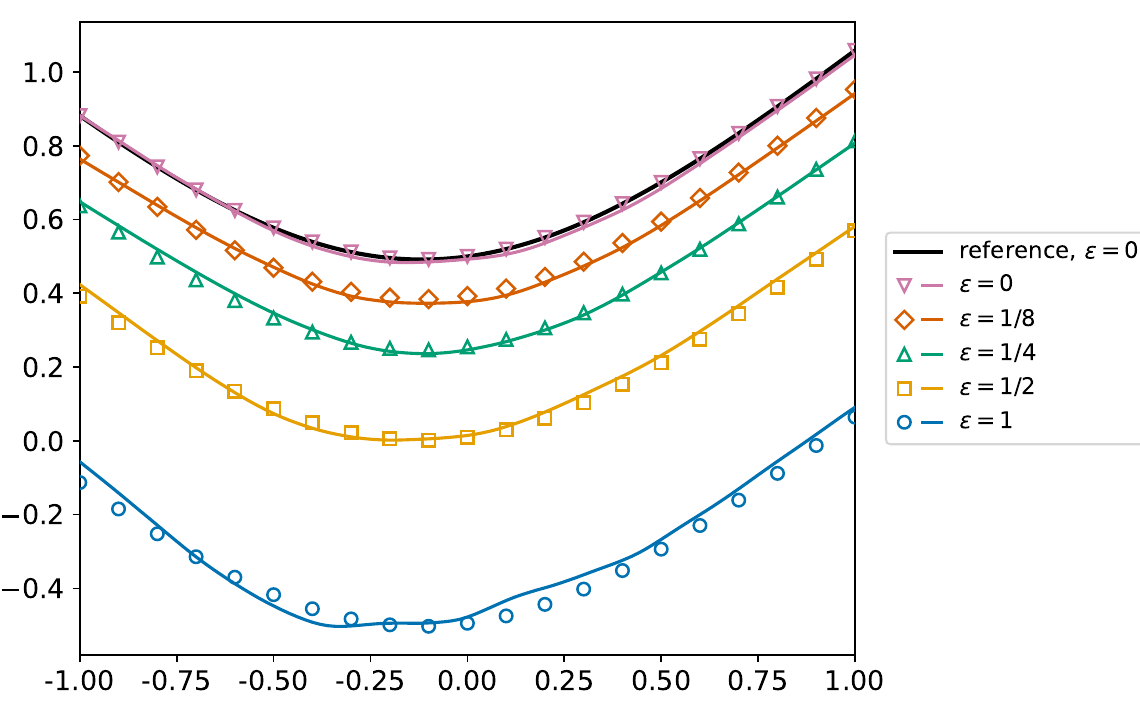}
    \caption{
        Profiles of the learned value function $v_{\theta}$ and of the empirical costs $\hat{J}(u_{\alpha})$ for the SOCP in \Cref{sec_eps}, with fixed $d = 10^4$ and varying perturbation parameter $\epsilon$.
    Lines depict $s \mapsto v_{\theta}(0, s\B{1}_d)$, whereas markers in the corresponding colors depict $\hat{J}(u_{\alpha})$ defined in \eqref{eq_empJ} with the initial condition $\bar{X}_0^m = s\B{1}_d$.    }\label{fig_socp_eps}
\end{figure}

As $\epsilon$ decreases, the solution of \eqref{eq_hjb_eps} is expected to converge to the reference solution at $\epsilon = 0$.
The numerical results are presented in \Cref{fig_socp_eps}.
As $\epsilon$ decreases, the learned value function approaches the reference solution for $\epsilon = 0$, as expected. For each fixed $\epsilon$, the learned value function closely matches the empirical cost over the displayed initial states, indicating that $v_{\theta}$ accurately represents the cost attained by the learned feedback control $u_{\alpha}$. 
These results demonstrate the effectiveness of the proposed method for SOCPs whose Hamiltonian 
infimization $\inf_u H$ does not admit an explicit closed-form expression.

\section{Conclusions}\label{sec_conclu}

In this paper, we propose a highly parallel, derivative-free martingale deep learning method for high-dimensional SOCPs through their associated HJB equations and the verification theory for the optimal control. 
The method reformulates the HJB equation into two coupled subproblems: a martingale problem for the value function and an integral optimality condition for the feedback control. 
As a result, the loss function avoids evaluating the gradient and Hessian of the value function. 
To enable parallel computation, we introduce a pilot process for state-space exploration and system processes for martingale enforcement, allowing the loss to be evaluated in parallel across both time and space. 
The martingale condition is imposed in a weak Galerkin form via adversarial learning, which further improves training efficiency and robustness. 
The optimal feedback control is learned from the integral optimality condition, thereby avoiding pointwise minimization for the optimal control.
As a result, an efficient training algorithm is proposed to solve the value function and the optimal control simultaneously for high-dimensional SOCPs. 

Numerical experiments on a range of high-dimensional problems demonstrate that the proposed method is both accurate and efficient in dimensions up to $10^4$, including cases in which the optimal control does not admit an explicit closed-form expression.
The results also indicate that network expressivity can become a practical bottleneck for extremely high-dimensional problems with complicated solution profiles. 
Future work will focus on the theoretical analysis of the proposed scheme, and the development of more efficient network architectures to further enhance the performance of the method for large-scale problems.

% \section*{Acknowledgments}
% We would like to acknowledge the assistance of volunteers in putting
% together this example manuscript and supplement.

% \section*{Code and data availability}
% The code and data that support the findings of this study are available from the corresponding author upon reasonable request.

\bibliographystyle{siamplain}
\bibliography{bibliography}

\end{document}

%% file: ex_shared.tex
\usepackage{lipsum}
\usepackage{amsfonts}
\usepackage{graphicx}
\usepackage{epstopdf}
\ifpdf
  \DeclareGraphicsExtensions{.eps,.pdf,.png,.jpg}
\else
  \DeclareGraphicsExtensions{.eps}
\fi

\newsiamremark{remark}{Remark}
\newsiamremark{hypothesis}{Hypothesis}
\crefname{hypothesis}{Hypothesis}{Hypotheses}
\newsiamthm{claim}{Claim}
\newsiamremark{fact}{Fact}
\crefname{fact}{Fact}{Facts}

\headers{Derivative Free Martingale Network}{W. Cai, S.X. Fang, W.Z. Zhang, T. Zhou}

\title{Martingale deep neural network for very high-dimensional stochastic optimal controls\thanks{Date: \today \;(First version available on arXiv:2408.14395 as preprint	on Aug. 26, 2024).
\funding{
The work of WC is funded by Clements Chair of Applied Math. at SMU.  
The work of SF is supported by the Postdoctoral Fellowship Program of CPSF (under grant GZB20250705) and China Postdoctoral Science Foundation (under grant 2025M783150). 
The work of WZ is supported by the National Key R\&D Program of China (under grant 2022YFA1005203) and the NSF of China (under grant 92270205, grant 1220012530).  
The work of TZ is supported by the NSF of China (under grant 12288201) and the Youth Innovation Promotion Association (CAS). 
}}}

\author{Wei Cai\thanks{Corresponding author. Department of Mathematics, Southern Methodist University, 208 Clements Hall, 3100 Dyer St, Dallas, TX 75205, USA. (\email{cai@smu.edu}).}
\and Shuixin Fang\thanks{Institute of Computational Mathematics and Scientific/Engineering Computing, Academy of Mathematics and Systems Science, Chinese Academy of Sciences, Beijing, China. (\email{sxfang@amss.ac.cn}, \email{tzhou@lsec.cc.ac.cn}).}
\and Wenzhong Zhang\thanks{Suzhou Institute for Advanced Research, University of Science and Technology of China, Suzhou, China (\email{wenzhong@ustc.edu.cn}).}
\and Tao Zhou\footnotemark[3]
}

\usepackage{amsopn}

\usepackage[sort]{cite}

\usepackage{booktabs}
\usepackage{multirow}

\usepackage{algpseudocode}
\algrenewcommand\algorithmicrequire{\textbf{Input:}}
\algrenewcommand\algorithmicensure{\textbf{Output:}}

\usepackage{float}
\usepackage[caption=false]{subfig}
\usepackage{pgfplots}
\pgfplotsset{compat=1.18}
\usepgfplotslibrary{groupplots}
\usepgfplotslibrary{fillbetween}

\newcommand{\R}{\mathbb{R}}
\newcommand{\di}{\,\mathrm{d}}

\newcommand{\M}{\mathcal{M}}

\newcommand{\mO}{\mathcal{O}}

\newcommand{\E}[1]{\mathbb{E}\left[#1\right]}

\newcommand{\br}[1]{\left(#1\right)}
\newcommand{\fbr}[1]{(#1)}
\newcommand{\abs}[1]{\left\vert#1\right\vert}

\newcommand{\bbr}[1]{\left\{#1\right\}}

\newcommand{\vbr}[2]{#1(#2 #1)}

\newcommand{\argmin}{\mathop{\arg\min}}

\newcommand{\B}[1]{\boldsymbol{#1}}
\newcommand{\sptext}[1]{\;\; \text{#1} \;\;}

%% file: fig_scaling_tikz.tex
    \begin{tikzpicture}
      \definecolor{cTime}{HTML}{0072B2}
      \definecolor{cMem}{HTML}{D55E00}
      \begin{axis}[
          name=A,
          width=6.4cm, height=5.4cm,
          xmode=log, ymode=log,
          xmin=60, xmax=1500,
          xtick={75,150,300,600,1200},
          xticklabels={75,150,300,600,1200},
          ymin=0.0052, ymax=0.0125,
          ytick={0.006,0.008,0.010},
          yticklabels={6E-3,8E-3,1E-2},
          xlabel={$r$}, ylabel={$\mathrm{RE}_{v}^{t_0}$},
          tick label style={font=\footnotesize}, label style={font=\small},
        ]
        \addplot[name path=hi, draw=none, forget plot] coordinates
          {(75,0.011487) (150,0.008479) (300,0.007712) (600,0.006990) (1200,0.006382)};
        \addplot[name path=lo, draw=none, forget plot] coordinates
          {(75,0.0095845) (150,0.0068053) (300,0.0066100) (600,0.0059805) (1200,0.0060498)};
        \addplot[cTime!25, forget plot] fill between[of=hi and lo];
        \addplot[cTime, mark=*, mark size=1.7pt, line width=0.9pt] coordinates
          {(75,0.0095845) (150,0.0068053) (300,0.0066100) (600,0.0059805) (1200,0.0060498)};
      \end{axis}
      \begin{axis}[
          name=B, at={(A.east)}, anchor=west, xshift=1.5cm,
          width=6.4cm, height=5.4cm,
          xmode=log,
          xmin=60, xmax=1500,
          xtick={75,150,300,600,1200},
          xticklabels={75,150,300,600,1200},
          ymin=0.994, ymax=1.10,
          ytick={1.00,1.05,1.10},
          xlabel={$r$}, ylabel={ratio to $r = 75$},
          tick label style={font=\footnotesize}, label style={font=\small},
          legend style={font=\scriptsize, draw=black!30, inner sep=1.5pt},
          legend pos=north west,
        ]
        \addplot[black!45, thin, forget plot] coordinates {(60,1) (1500,1)};
        \addplot[cTime, mark=*, mark size=1.7pt, line width=0.9pt] coordinates
          {(75,1.0000) (150,1.0025) (300,1.0023) (600,1.0069) (1200,1.0130)};
        \addlegendentry{RT}
        \addplot[cMem, mark=square*, mark size=1.6pt, line width=0.9pt] coordinates
          {(75,1.0000) (150,1.0005) (300,1.0014) (600,1.0033) (1200,1.0072)};
        \addlegendentry{GPU mem.}
      \end{axis}
    \end{tikzpicture}

%% file: fig_cost_tikz.tex
    \begin{tikzpicture}
      \definecolor{cTime}{HTML}{0072B2}
      \definecolor{cMem}{HTML}{D55E00}
      \begin{groupplot}[
          group style={group size=4 by 1, horizontal sep=0pt},
          scale only axis,
          width=2.80cm, height=3.0cm,
          xmode=log, ymode=log,
          log ticks with fixed point,
          every axis plot/.append style={line width=0.9pt, mark size=1.7pt},
          tick label style={font=\footnotesize},
          label style={font=\small},
        ]
        \nextgroupplot[
          xlabel={$N$}, ylabel={Relative cost},
                    xmin=7.2, xmax=220,
          xtick={8,16,32,64,128},
                    xticklabels={$2^3$,$2^4$,$2^5$,$2^6$,$2^7$},
          ymin=0.93, ymax=3.5, ytick={1,2,3},
          ymajorgrids,
          major grid style={black!45, dashed, thin},
        ]
        \addplot[cTime, mark=*] coordinates {
          (8,1.0000) (16,1.0043) (32,1.0296) (64,1.0380) (128,1.0586) };
        \addplot[cMem, mark=square*] coordinates {
          (8,1.0000) (16,1.0459) (32,1.1377) (64,1.5393) (128,2.6421) };
        \node[draw=black!30, inner sep=0.8pt, font=\tiny, anchor=north west, align=left,
    xshift=-32pt]
    at (axis cs:7.6,3.00) {%
            \tikz{\draw[cTime, line width=0.7pt] plot[mark=*, mark size=1.3pt] coordinates {(0,0) (0.42,0)};}\,RT per iter.\\[0.8pt]
            \tikz{\draw[cMem, line width=0.7pt] plot[mark=square*, mark size=1.2pt] coordinates {(0,0) (0.42,0)};}\,GPU mem.};
        \nextgroupplot[
                    xlabel={$\abs{\mathbb{A}_1}$, $\abs{\mathbb{A}_2}$}, xmin=20, xmax=900,
          xtick={32,64,128,256,512},
                    xticklabels={$2^5$,$2^6$,$2^7$,$2^8$,$2^9$},
          ymin=0.93, ymax=3.5, ytick={1,2,3},
          yticklabels={,,},
          ymajorgrids,
          major grid style={black!45, dashed, thin},
        ]
        \addplot[cTime, mark=*] coordinates {
          (32,1.0000) (64,1.1188) (128,1.3637) (256,1.9693) (512,3.1669) };
        \addplot[cMem, mark=square*] coordinates {
          (32,1.0000) (64,1.0013) (128,1.0040) (256,1.0094) (512,1.0198) };
    \nextgroupplot[xlabel={$M$}, xmin=600, xmax=18200,
    xtick={1024,2048,4096,8192,16384},
                            xticklabels={$2^{10}$,$2^{11}$,$2^{12}$,$2^{13}$,$2^{14}$},
    ymin=0.93, ymax=3.5, ytick={1,2,3},
    yticklabels={,,},
    ymajorgrids,
    major grid style={black!45, dashed, thin},
    ]
    \addplot[cTime, mark=*] coordinates {
        (1024,1.0000) (2048,1.0005) (4096,1.0008) (8192,1.0002) (16384,1.0009) };
    \addplot[cMem, mark=square*] coordinates {
        (1024,1.0000) (2048,1.0194) (4096,1.0580) (8192,1.1352) (16384,1.4195) };
    \nextgroupplot[
                xshift=0.70cm,
                xlabel={$d$}, xmin=930, xmax=18200,
                xtick={1024,2048,4096,8192,16384},
                xticklabels={$2^{10}$,$2^{11}$,$2^{12}$,$2^{13}$,$2^{14}$},
                ymin=0.9, ymax=60, ytick={1,2,5,10,20,50},
                extra y ticks={1},
                extra y tick style={grid=major, grid style={black!45, thin}, yticklabel={}},
    ]
    \addplot[cTime, mark=*] coordinates {
                (1024,1.0000) (2048,1.3817) (4096,2.7768) (8192,9.5737) (16384,33.8884) };
    \addplot[cMem, mark=square*] coordinates {
                (1024,1.0000) (2048,2.1258) (4096,4.7854) (8192,12.6670) (16384,44.2965) };
      \end{groupplot}
    \end{tikzpicture}